\documentclass{article}[15pt]
\usepackage{amsmath, epsfig,amssymb, multicol}
\usepackage[active]{srcltx}
\usepackage{epstopdf}
\usepackage{color}
\newtheorem{lemma}{Lemma}
\newtheorem{theorem}{Theorem}
\newtheorem{conjecture}{Conjecture}

\newtheorem{claim}{Claim}
\newtheorem{fact}{Fact}

\newtheorem{remark}{Remark}

\title{The fractional chromatic number of $K_{\Delta}$-free graphs}
\author{Xiaolan Hu\thanks{School of Mathematics and Statistics \& Hubei Key Laboratory of Mathematical Sciences, Central China Normal University, Wuhan 430079, P.~R.~China. E-mail: {\tt xlhu@mail.ccnu.edu.cn},  Research is partially supported
by National Natural Science Foundation of China (No. 11971196).}
\and Xing Peng\thanks{School of Mathematical Sciences, Anhui University,  Hefei 230601, P.~R.~China. E-mail: {\tt x2peng@ahu.edu.cn}, Research is supported in part by National Natural Science Foundation of China (No.\ 12071002) and a Start--up Fund from Anhui University.}
}
\date{}

\begin{document}
\maketitle

\begin{abstract}
For a simple graph $G$, let $\chi_f(G)$ be the fractional chromatic number of $G$.
In this paper, we aim to establish upper bounds on $\chi_f(G)$ for those graphs $G$ with restrictions on the clique number. Namely, we prove that for $\Delta \geq 4$, if  $G$ has maximum degree at most $\Delta$ and is $K_{\Delta}$-free,  then $\chi_f(G) \leq \Delta-\tfrac{1}{8}$ unless $G= C^2_8$ or  $G = C_5\boxtimes K_2$. This  improves the result in   [King, Lu, and Peng, {\it SIAM J.~Discrete Math.,} 26(2) (2012), pp. 452-471] for $\Delta \geq 4$ and the result  in [Katherine and King, {\it SIAM J.~Discrete Math.,} 27(2) (2013), pp. 1184-1208] for $\Delta \in \{6,7,8\}$.
\end{abstract}


\section{Introduction}
\begin{figure}[htbp]
 \centerline{\psfig{figure=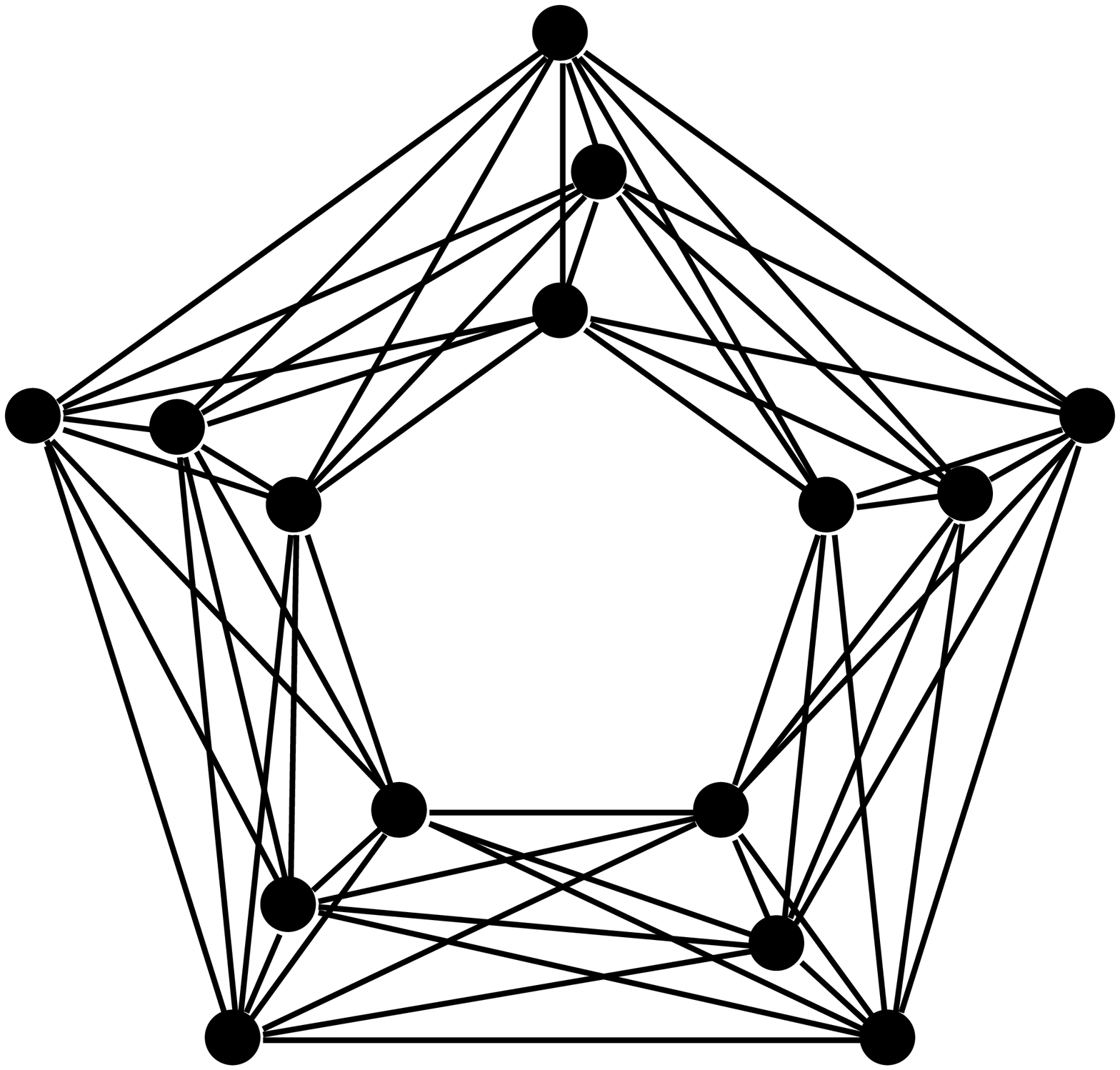, width=0.2\textwidth} \hfil
  \psfig{figure=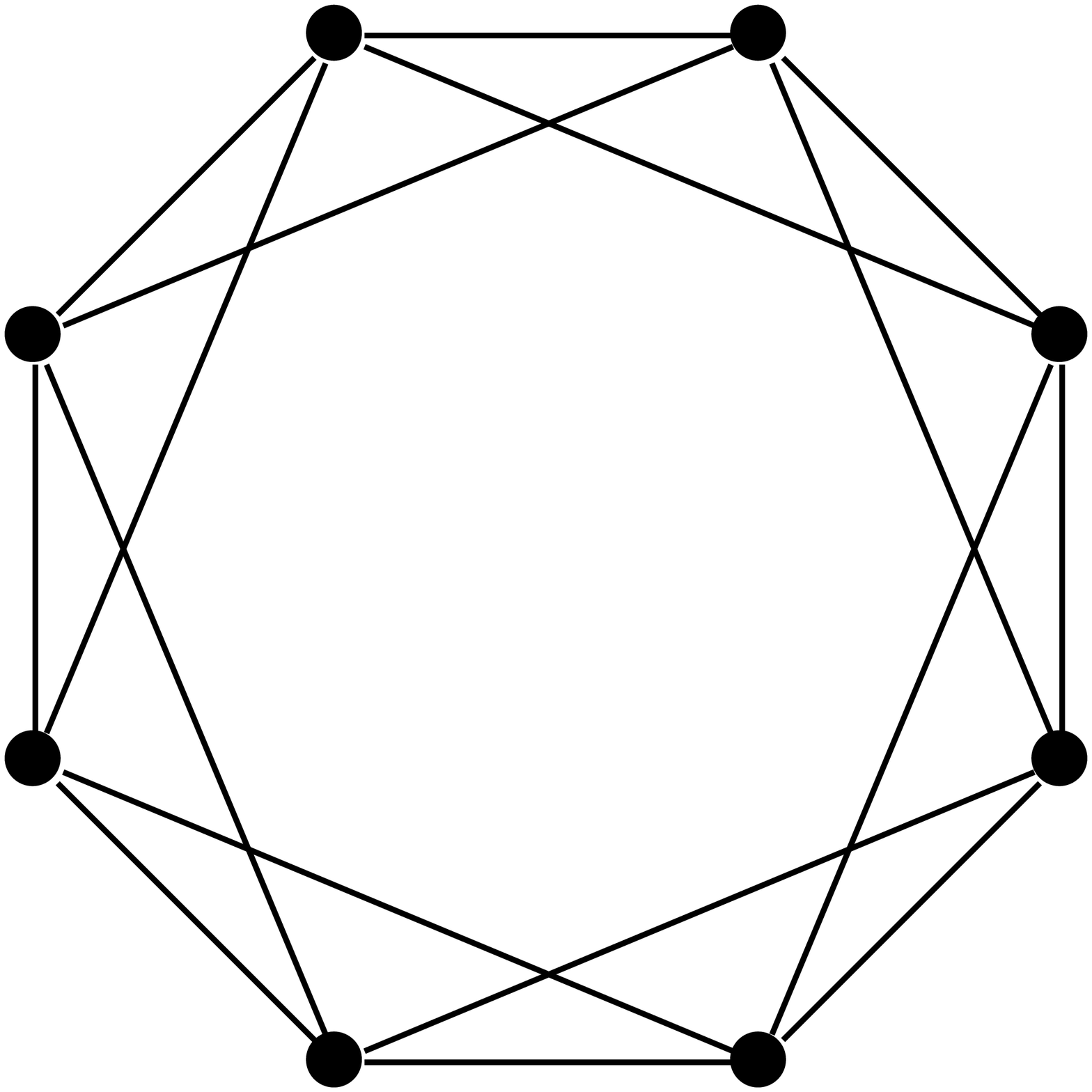, width=0.2\textwidth} \hfil
   \psfig{figure=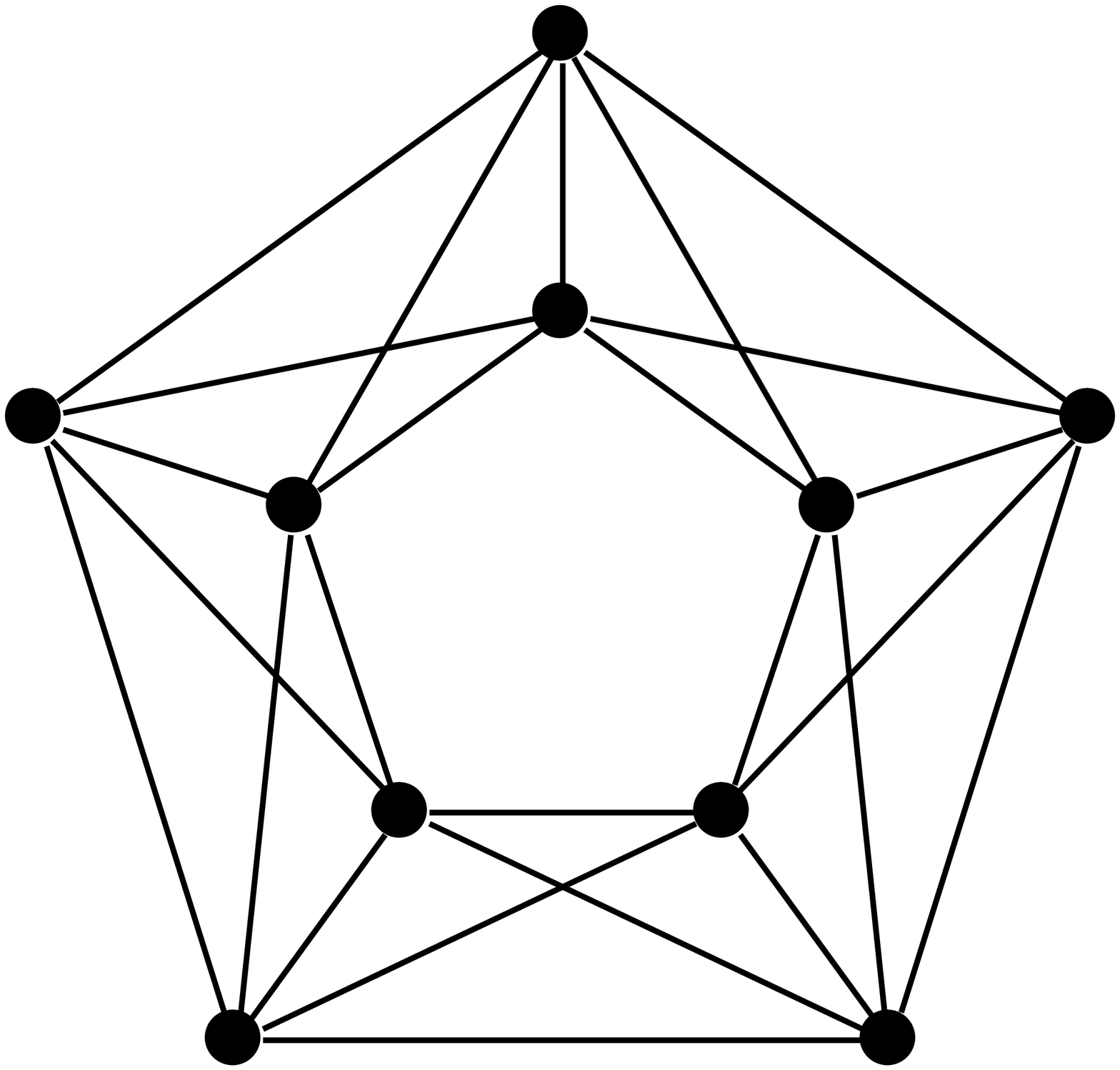, width=0.2\textwidth}}
 \centerline{\hspace{3.5cm}$C_5 \boxtimes K_3$ \hspace{3cm} $C_8^2$  \hspace{3.5cm} $C_5 \boxtimes K_2$\hspace{3.5cm}}
\caption{The graphs $C_5 \boxtimes K_3$, $C_8^2$, and $C_5 \boxtimes K_2$.}
 \label{examples}
\end{figure}

In this paper, we use $\chi(G)$ to denote the chromatic number of a simple graph $G$. It is well-known that $\Delta+1$ is a trivial upper bound for $\chi(G)$, where $\Delta$ is the maximum degree of $G$.  The famous Brooks' Theorem states that if the maximum degree of a connected graph $G$ is $\Delta$, then $\chi(G) \leq \Delta$ unless $G$ is $K_{\Delta+1}$ or an odd cycle. This seminal result leads to the idea bounding the chromatic number $\chi(G)$ in terms of the maximum degree $\Delta(G)$ and the clique number $\omega(G)$.
For $\Delta \geq 3$, Brooks' Theorem asserts that if $G$ is $K_{\Delta+1}$-free (i.e.,~$\omega(G) \leq \Delta(G)$), then $\chi(G) \leq \Delta(G)$. One may ask to show an upper bound for $\chi(G)$ if $G$ is $K_{\Delta}$-free. This is exactly the following conjecture by Borodin and Kostochka \cite{bk}.
\begin{conjecture}\label{bkc}
For a connected graph $G$,  if $\omega(G)\leq \Delta(G)-1$
and $\Delta(G)\geq 9$, then we have $\chi(G)\leq \Delta(G)-1.$
\end{conjecture}
If the conjecture is true, then it is best possible since there is a
$K_8$-free graph $G=C_5\boxtimes K_3$ (see Figure \ref{examples}) with $\Delta(G)=8$, $\omega(G)=6$, and $\chi(G)=8$.

Given two graphs $G$ and $H$, the {\it strong product} $G\boxtimes H$ is the graph with vertex set $V(G)\times V(H)$, and  $(u,x)$ is connected to $(v,y)$ if one of the following holds
\begin{description}
\item{1:} $u=v$ and $xy \in E(H)$,
\item{2:} $uv \in E(G)$ and $x=y$,
\item{3:} $uv \in E(G)$ and $xy\in E(H)$.
\end{description}

Using the probabilistic methods, Reed \cite{reed1} confirmed Conjecture \ref{bkc} for $\Delta \geq 10^{14}$. In the same paper, Reed claimed that a careful calculation could reduce $10^{14}$ to $10^3$. Let us assume that Conjecture \ref{bkc} is true for $\Delta \geq 10^3$. This is the state of the art of Conjecture \ref{bkc}. King, Lu, and Peng first considered to bound the fractional chromatic number of $K_{\Delta}$-free graphs, which can be viewed as the fractional analog of  Conjecture \ref{bkc}.

The fractional chromatic number $\chi_f(G)$ of a graph $G$ can be defined as follows. For positive integers $a$ and $b$ with $ b \leq a$, we use $[a]$ to denote the set of first $a$ positive integers and $\binom{[a]}{b}$ to denote the set of $b$-subsets of $[a]$.
An ($a$:$b$)-{\it
coloring} of a graph $G$ is a mapping $\psi: V(G) \rightarrow \binom{[a]}{b}$ such that
$\psi(u) \cap \psi(v)=\emptyset$ whenever $uv \in E(G)$.
The $b$-{\it fold coloring number} of $G$,
denoted by $\chi_b(G)$, is the smallest integer $a$ such that $G$
has an ($a$:$b$)-coloring. Note that $\chi_1(G)=\chi(G)$.   The {\it
  fractional chromatic number} $\chi_f(G)$ of $G$ is
 $\lim_{b \rightarrow \infty} \frac{\chi_b(G)}{b}.$  If $G$ is vertex transitive, then $\chi_f(G)=\tfrac{|V(G)|}{\alpha(G)}$, where $\alpha(G)$ is the independence number of $G$.

For $\Delta \geq 3$, we define
\[
f(\Delta)=\underset{G}\min \{\Delta-\chi_f(G): \Delta(G) \leq \Delta, ~\omega(G) \leq \Delta-1, ~ G \notin \{C_8^2, C_5 \boxtimes K_2\}\}.
\]

The value of $f(3)$ received a great deal of attention. Actually, based on the study of independent sets in triangle-free subcubic graphs \cite{faj,gm,ht,staton}, Heckman and Thomas \cite{ht} conjectured $f(3)=\tfrac{1}{5}$. Before  Dvo\v r\'ak, Sereni, and Volec \cite{DSV}  proved Heckman and Thomas' Conjecture, progress towards this conjecture can be found in
 \cite{FKK,hz,KKV,chl,LP}.
King, Lu, and Peng \cite{KLP} first proved that $f(\Delta) \geq \tfrac{2}{67}$ for each $\Delta \geq 4$,. Edwards and King \cite{EK} established an improved  lower bound on $f(\Delta)$ for all $\Delta \geq 6$. In particular, they showed $f(6) \geq \tfrac{1}{22.5}$, $f(7) \geq \tfrac{1}{11.2}$, $f(8) \geq \tfrac{1}{8.9}$, and $f(\Delta) > \tfrac{1}{7.7}$ for each $\Delta \geq 9$.

We prove the following theorem in this paper.
\begin{theorem}\label{main}
For each $\Delta \geq 4$, if $G$ has maximum degree at most $\Delta$ and is $K_{\Delta}$-free, then $\chi_f(G) \leq \Delta-\tfrac{1}{8}$ unless  $G= C^2_8$ or  $G = C_5\boxtimes K_2$ (see Figure \ref{examples}).
\end{theorem}
 By the theorem above, we get $f(\Delta) \geq \tfrac{1}{8}$ for each $\Delta \geq 4$. Therefore, Theorem \ref{main} improves the result of King, Lu, and Peng \cite{KLP} for each $\Delta \geq 4$ and the result of Edwards and King \cite{EK} for $\Delta\in \{6,7,8\}$.

The rest of this paper is organized as follows. In Section 2, we will prove Theorem \ref{main}. Section 3 is devoted to prove technical lemmas.

\section{Proof of Theorem \ref{main}}
Recall the definition of $f(\Delta)$.  King, Lu, and Peng \cite{KLP} proved the following lemma on $f(\Delta)$.
 \begin{lemma}\label{increase}
For each $\Delta \geq 6$, we have $f(\Delta)\geq \min \left\{f(\Delta-1),
\frac{1}{2}\right\}.$ We also have  $f(5)\geq \min\left\{f(4),
\frac{1}{3}\right\}$.
\end{lemma}
 In view of Lemma \ref{increase}, in order to prove Theorem \ref{main}, we need only to show $f(4) \geq \tfrac{1}{8}$. We define
 \[
 {\cal G}=\{G: \Delta(G) \leq 4,~ \omega(G) \leq 3, \textrm{ and } G \not =C_8^2\}.
 \]
  Equivalently, we have to show $\chi_f(G) \leq 4-\tfrac{1}{8}$ for each $G \in {\cal G}$.
Actually, we will prove a stronger result. To state our result, we need to introduce some definitions.

Let $G$ be a graph. For $v\in V(G)$, the {\it degree} $d_G(v)$ of $v$ is the number of edges incident with $v$ in $G$.
For two vertices $u$ and $v$ of $G$, let $\textrm{dist}_G(u,v)$ be the length of a shortest path between $u$ and $v$.
For $v \in V(G)$ and a positive integer $j $,
let $N_G^j(v)=\{u \in V(G): \textrm{dist}_G(u,v)=j\}$. If $S$ and $T$ are two  disjoint subsets of $V(G)$, then the number of edges between $S$ and $T$ is denoted by $e(S,T)$.
A {\it demand function} $h_G$ is a function from $V(G)$ to $[0,1]$ taking rational values. For  positive integers $a$ and $b$ with $a <b$, we use $[a,b]$ to denote the set $\{a,a+1,\ldots,b\}$. If $a=1$, then we write $[b]$ instead of $[1,b]$ for short. Suppose $h_G$ is a demand function and  $N$ is a positive integer such that $N\cdot h_G(v)$ is a positive integer for each $v \in V(G)$. Let $\psi$ be a function from $V(G)$ to the set of all subsets of $[N]$. If $|\psi(v)| \geq Nh_G(v)$ holds for each $v \in V(G)$, then we say $\psi$ is an $(h_G,N)$-coloring of $G$. We mention a few remarks on these definitions.
\begin{description}
\item (i) If $G$ has an $(h_G,N)$-coloring, then $G$ also has an $(h_G,M)$-coloring provided $M$ is divisible by $N$. Actually, we can duplicate each color from an $(h_G,N)$-coloring $M/N$ times to get an $(h_G,M)$-coloring.
\item (ii) If $\psi$ is an  $(h_G,N)$-coloring of $G$, then we can assume that $|\psi(v)|=Nh_G(v)$ for each $v \in V(G)$ by deleting colors from $\psi(v)$. We will keep this assumption in the rest of this paper.
\item (iii) To show $\chi_f(G) \leq r$ for a rational number $r$, it is enough to show $G$ has an $(h_G,N)$-coloring such that $h_G(v) \geq 1/r$ for each $v \in V(G)$.
\end{description}
In order to prove $\chi_f(G) \leq 4-\tfrac{1}{8}$ for each $G \in {\cal G}$, we define a demand function $h_G$ as following:
\begin{description}
\item (i) If $1 \leq d_G(v) \leq 4$, then  $h_G(v)=\tfrac{12-d_G(v)}{31}$;
\item (ii) If $d_G(v)=0$, then  $h_G(v)=1$.
\end{description}
We will prove the following stronger result.
\begin{theorem}\label{stronger}
For each $G \in {\cal G}$, there is a positive integer $N_G$ such that $G$ has an $(h_G,N_G)$-coloring.
\end{theorem}
It is easy to see Theorem \ref{stronger} implies that $\chi_f(G) \leq 4-\tfrac{1}{8}$  for each $G \in {\cal G}$. We next prove Theorem \ref{stronger} by contradiction. Let $G \in \cal G$ be a counterexample with the smallest number of vertices.
The following simple remark follows from the remark (i) on the definition of $(h_G,N)$-coloring.
\noindent
\begin{remark}
 There exists a positive integer $M$ such that each graph $G'$ with $|V(G')|<|V(G)|$ has an $(h_{G'},M)$-coloring.
\end{remark}
As $G$ is a minimum counterexample, $G'$ has an $(h_{G'},N_{G'})$-coloring for each $G'$ with $|V(G')|<|V(G)$. Thus we can take $M$ as the least common multiplier of $N_{G'}$, where $G'$ ranges over all graphs  satisfying $|V(G')|<|V(G)|$.

 We will show that $G$ satisfies the following properties.
\begin{lemma} \label{key0}
The graph $G$ is $2$-connected and  satisfies $\delta(G)  \geq 3$.
\end{lemma}
\noindent
{\bf Proof:} Clearly, $G$ is connected. If $G$ contains a cut-vertex $v$, then there are two subgraphs $G_1$ and $G_2$ such that
$V(G_1)\cup V(G_2)=V(G)$ and $V(G_1) \cap V(G_2)=\{v\}$. As $G$ is a minimum counterexample,  the graph $G_i$ has  an $(h_{G_i},M)$-coloring $\psi_i$ for $1 \leq i \leq 2$. By deleting colors if it is necessary, we can assume that $|\psi_1(v)|=|\psi_2(v)|=M \cdot h_G(v)$, here $h_G(v)=\tfrac{12-d_G(v)}{31}$. By permuting colors in $\psi_2$, we can further assume that $\psi_1(v)=\psi_2(v)$. Thus the union of $\psi_1$ and $\psi_2$ is an $(h_G,M)$-coloring of $G$. A contradiction.

Since $G$ is 2-connected, we get $\delta(G)\geq 2$. Suppose that there exists $v\in V(G)$ with $d_G(v)=2$. Let $u$ and $w$ be the neighbors of $v$. By the minimality of $G$, the graph $G'=G-v$ has an $(h_{G'},M)$-coloring, say $\psi$.
Deleting colors from $\psi(u)$ and $\psi(w)$, we can assume that $|\psi(u)|=M \cdot \tfrac{12-d_G(u)}{31}$ and $|\psi(w)|=M \cdot \tfrac{12-d_G(w)}{31}$. Since $d_G(u)\geq 2$ and $d_G(w)\geq 2$, it follows that $|\psi(u)|\leq \tfrac{10}{31} \cdot M$ and $|\psi(w)|\leq \tfrac{10}{31} \cdot M$. As $|[M] \setminus (\psi(u) \cup \psi(w))| \geq \tfrac{11}{31} \cdot M \geq h_G(v) \cdot M$, we can extend the coloring $\psi$ to $v$ and obtain an $(h_G,M)$-coloring of $G$.  A contradiction. \hfill $\square$
%
%

~~~

We also need the following lemma on the local structures of $G$.
\begin{lemma} \label{key4}
For each $v \in V(G)$, if $d_G(v)=3$, then we have $e(N_G(v),N_G^2(v)) \geq 5$;  if  $d_G(v)=4$, then we have $e(N_G(v),N_G^2(v)) \geq 9$.
\end{lemma}
The proof of Lemma \ref{key4} is quite involved, we postpone it until Section 3 and proceed to prove Theorem \ref{stronger}.

~~~

\noindent
{\bf Proof of Theorem \ref{stronger}:} We will prove that $G$ has an $(h_G,M')$-coloring for some positive integer $M'$.  This will be a  contradiction to the assumption of $G$. For each $v \in V(G)$, let $H_v$ be the proper subgraph of $G$ obtained from $G$ by deleting vertices from $N_G(v)$. By Remark 1 and the choice of $G$, the graph $H_v$ has an $(h_{H_v},M)$-coloring, say $\psi_v$. Let $\psi$ be the union of  $\psi_v$, where $v$ ranges over all vertices of $G$ and colors used for different $\psi_v$ are distinct. Let $n=|V(G)|$. We next show $\psi$ is an $(h_G,M')$-coloring, where $M'=nM$. To show this, we have to verify $|\psi(u)| \geq M'h_G(u)$ for each vertex $u$.  By Lemma \ref{key0}, there are two cases depending on the degree of $u$.

If $d_G(u)=3$, then we write $N_G^2(u)=\cup_{i=1}^3 A_i$, where $w \in A_i$ if $w$ has $i$ neighbors in $N_G(u)$.
We can compute the number of colors $u$ received as follows. In the graph $H_u$, the vertex $u$ is an isolated one and it receives $M$ colors from $\psi_u$. For each $v \in N_G(u)$, the vertex $u$ is not in the graph $H_v$ and does not receive any color from $\psi_v$. If $v \in N_G^2(u) \cap A_1$, then $u$ has degree two in $H_v$ and $u$ receives $10M/31$ colors from $\psi_v$. Similarly, for $v \in N_G^2(u) \cap A_2$, the coloring $\psi_v$ gives $u$ exactly $11M/31$ color as $u$ has degree one in $H_v$. For $v \in N_G^2(u) \cap A_3$, the vertex $u$ receives $M$ colors from $\psi_v$ as $u$ is an isolated vertex in $H_v$. If $v \in V(G)\setminus (\{u\} \cup N_G(u) \cup N_G^2(u))$, then $u$ has degree three in $H_v$ and so it receives $9M/31$ colors from $\psi_v$.
For $1 \leq i \leq 3$, let $a_i=|A_i|$.
The number of colors for $u$ can be calculated as follows:
\begin{align*}
|\psi_u|&=M+a_1 \cdot (10M)/31+a_2\cdot (11M)/31+a_3\cdot M\\
         &\hspace{0.4cm}  +(n-4-a_1-a_2-a_3)\cdot (9M)/31\\
       &= (9n+a_1+2a_2+22a_3-5) \cdot M/31 \\
       & \geq (9n+a_1+2a_2+3a_3-5) \cdot M/31 \\
       & \geq (9nM)/31,
\end{align*}
in the last step, we notice   $e(N_G(u),N_G^2(u))=a_1+2a_2+3a_3 \geq 5$ as the first part of Lemma \ref{key4}.

If $d_G(u)=4$, then we write $N_G^2(u)=\cup_{i=1}^4 A_i$, where $w \in A_i$ if $w$ has $i$ neighbors in $N_G(u)$.  Let $a_i=|A_i|$ for each $1 \leq i \leq 4$. As we did in the previous case, we can estimate the size of $\psi_u$ as follows:
\begin{align*}
|\psi_u|&=M+a_1 \cdot (9M)/31+a_2\cdot (10M)/31+a_3\cdot (11M)/31+a_4 \cdot M\\
         &\hspace{0.4cm}  +(n-5-a_1-a_2-a_3-a_4)\cdot (8M)/31\\
         &= (8n+a_1+2a_2+3a_3+23a_4-9) \cdot M/31 \\
         & \geq (8n+a_1+2a_2+3a_3+4a_4-9) \cdot M/31 \\
         & \geq (8nM)/31,
\end{align*}
in the last step, we note that $e(N_G(u),N_G^2(u))=a_1+2a_2+3a_3+4a_4 \geq 9$ by the second part of Lemma \ref{key4}.

Therefore, $\psi$ is an $(h_G,M')$-coloring for $G$. This is a contradiction and  the theorem is proved. \hfill $\square$

\section{Proof of Lemma \ref{key4}}

We start with some definitions and notations. Let $G$ be a graph.  A subset $S\subseteq V(G)$ is called an {\it independent set} if there are no edges with both endpoints in $S$.
The {\it contraction} of  $S$ is to replace vertices in
$S$ by a single fat vertex, denoted by $\underline{S}$, whose incident edges are all edges that are incident with at least one vertex in $S$, except edges with both endpoints in $S$. The new
graph obtained by contracting $S$ is denoted by $G/S$. This operation is also known as
identifying vertices of $S$ in the literature. For convenience, if $S=\{u,v\}$, then the fat vertex will be denoted by $\underline{uv}$.
We use $G-S$ to denote the subgraph of $G$ induced by $V(G)-S$.
For two nonadjacent vertices $u$ and $v$, let $G+uv$ be the graph obtained from $G$ by adding an edge $uv$.
The square $G^2$ of a graph $G$ is the graph with vertex set $V(G)$ and two vertices are adjacent in $G^2$  if and only if their distance in $G$ is at most $2$.

Recall that $G$ is a  counterexample with the smallest number of vertices.
 If $H$ is a subgraph of $G$ with $|V(H)|<|V(G)|$, then $H$ has an $(h_H,M)$-coloring $\psi$ by the definition of $G$.
 For each $v \in V(G)-V(H)$, we define $L(v)=M\setminus \cup_{u \in N_H(v)} \psi(u)$. Since $d_{H}(v)\leq d_{G}(v)$ for any $v\in V(H)$, we get $|\psi(v)|=(12-d_H(v))/31 \cdot M \geq (12-d_G(v))/31\cdot M$. If
$\psi'$ is an $(h_G,M)$ coloring of $H$ such that $\psi'(v)\subseteq \psi(v)$  for each $v \in V(H)$, then we say $\psi'$ is a  {\em limitation} of $\psi$ on $H$. Given a limitation $\psi'$ and a vertex $v \in V(G)-V(H)$, let $L'(v)=M\setminus \cup_{u \in N_H(v)} \psi'(u)$. Roughly speaking, if we have the partial coloring $\psi'$ of $G$ and we wish to color $u$, then $L'(u)$ is the set of colors which could be assigned to the vertex $u$. We also have the following simple fact which will be frequently used in our proof. Namely, we are able to prescribe a set of colors for $L'(u)$ by defining the limitation carefully.

\begin{fact} \label{fact1}
Let $S \subset V(G)$ and $H$ be the subgraph of $G$ induced by the vertex set $V(G)-S$. Suppose that $\psi$ is an $(h_H,M)$-coloring of $H$.
For $v_1,v_2 \in S$, if $A_1,A_2 \subset [1,M]$ with $|A_1| \leq M/31$ and $|A_2|\leq M/31$, then
there is a limitation $\psi'$ of $\psi$ on $H$ such that $A_1 \subset L'(v_1)$ and $A_2 \subset L'(v_2)$. Here, $A_1$ and $A_2$ are disjoint if $v_1$ and $v_2$ are adjacent.
\end{fact}
\noindent
{\bf Proof:} For each $u \in H$, if $N_G(u) \cap \{v_1,v_2\}=\emptyset$, then to define $\psi'(u)$, we delete colors from $\psi(u)$ arbitrarily if it is necessary.  If $|N_G(u) \cap \{v_1,v_2\}|=1$, say $u$ is adjacent to $v_1$, then we note that $h_H(u) \geq h_G(v)+1/31$. In order to define $\psi'(u)$, we need to delete  at least $M/31$ colors from $\psi(u)$. Therefore, we can  first remove colors of $A_1$ from $\psi(u)$  and then delete extra colors if it is necessary. If $|N_G(u) \cap \{v_1,v_2\}|=2$, then we have $h_H(u) \geq h_G(u)+2/31$. To define $\psi'(u)$, we remove colors of $A_1 \cup A_2$ from $\psi(u)$ first and  delete extra colors if it is necessary. It is easy to see that the limitation defined above  is a desired one.
\hfill\rule{1ex}{1ex}

 We remark that  one of $A_1$ and $A_2$ could be the empty set when we apply Fact \ref{fact1}.
Next, we show properties that satisfied by the minimal counterexample $G$.
\begin{fact}\label{fact2}
Suppose that $G$ contains a subgraph consisting of a triangle $uvw$ and two edges $uu',vv'$. Let $\psi$ be an $(h_{G'},M)$-coloring of $G'=G-\{u,v,w\}$. If $d_G(u)=d_G(v)=3$ and $|\psi(u')\setminus \psi(v')|\geq (8-d_G(w))M/31$, then any limitation of $\psi$ on $G'$ can be extended as an $(h_{G},M)$-coloring of $G$.
\end{fact}
\noindent
{\bf Proof:}  We choose $A \subset \psi(u')\setminus \psi(v')$ with $|A|=(8-d_G(w))M/31$. Let $\psi'$ be any limitation of $\psi$ on $G'$. Then $A\subseteq L'(v)$.
Note that $|L'(w)|\geq M-(d_G(w)-2)\cdot 9M/31$ and $|L'(w)\setminus A|\geq (41-8d_G(w))/31\geq (12-d_G(w))/31$.
We can first extend $\psi'$ to $w$ such that $\psi'(w)\cap A=\emptyset$, then extend  $\psi'$ to $v$ such that $A\subseteq \psi'(v)$. Notice that $|L'(u)|\geq M-(12-d_G(w))M/31-9M/31\times 2+(8-d_G(w))M/31=9M/31$. Thus we can color $u$ properly and  get an $(h_{G},M)$-coloring of $G$, a contradiction.  \hfill\rule{1ex}{1ex}

\begin{claim}
\label{cl:c82-e}
 $G+uv$ contains no $C_8^2$ for any $uv\notin E(G)$.
\end{claim}
\noindent {\bf Proof:} Suppose  that there exists $uv\notin E(G)$ such that
$G+uv$ contains a $C_8^2$.  Then $uv$ must be contained in this $C_8^2$.
Assume that $V(C_8^2)=\{v_1,v_2,\ldots,v_8\}$, where $v_iv_{i+1}$ and $v_iv_{i+2}$ are edges for each $1 \leq i \leq 8$ while the addition is under module 8.
Let $G'$ be the graph obtained from $G-\{v_1,v_2,\ldots,v_8\}\backslash \{u,v\}$ by identifying $u$ and $v$.
 Note that $\underline{uv}$ has degree at most two in $G'$. Therefore, $G'$ is $K_4$-free and $C_8^2$-free.
  Thus $G'$ has an $(h_{G'},M)$-coloring $\psi$ by the definition of $G$ and the fact $|V(G')|<|V(G)|$. Let $\psi'$ be a limitation of $\psi$ on $G-\{v_1,v_2,\ldots,v_8\}\backslash \{u,v\}$. Notice that $\min \{d_G(u),d_G(v)\} \geq 3$ and $d_{G'}(\underline{uv}) \leq 2$. We can suppose that $\psi'(u),\psi'(v)\subseteq [1,9M/31]$.
 Without loss of generality, we assume that $u=v_1$. Then either $v=v_2$ or $v=v_3$ by the symmetry.

 If $v=v_2$, then we extend $\psi'$ as follows.
  Let $\psi'(v_3)=[1+9M/31,17M/31]$, $\psi'(v_4)=[1+17M/31,25M/31]$, $\psi'(v_5)=[1,2M/31] \cup [1+25M/31,M]$, $\psi'(v_6)=[1+3M/31,11M/31]$, $\psi'(v_7)=[1+11M/31,19M/31]$, and $\psi'(v_8)=[1+19M/31,27M/31]$. We obtained an  $(h_G,M)$-coloring of $G$, which is a contradiction.

 If $v=v_3$, then we extend $\psi'$ as follows. Let $\psi'(v_2)=[1+23M/31,M]$, $\psi'(v_4)=[1+9M/31,17M/31]$, $\psi'(v_5)=[1+17M/31,25M/31]$, $\psi'(v_6)=[1,8M/31]$, $\psi'(v_7)=[1+26M/31,M] \cup [1+9M/31,12M/31]$, and $\psi'(v_8)=[1+14M/31,22M/31]$. Thus $G$ has an $(h_G,M)$-coloring of $G$, which is a contradiction.  \hfill\rule{1ex}{1ex}

\begin{claim}
\label{cl:c82-v}
 $G/\{x,y\}$ contains no $C_8^2$ for any  $xy \not \in E(G)$.
\end{claim}

\noindent {\bf Proof:} Suppose that $G/\{x,y\}$ contains a $C_8^2$ for some $xy \not \in E(G)$.
Then $\underline{xy}$ must be contained in this $C_8^2$.
Assume that $V(C_8^2)=\{v_1,v_2,\ldots,v_8\}$, where $v_iv_{i+1}$ and $v_iv_{i+2}$ are edges for each $1 \leq i \leq 8$ while the addition is under module 8. Without loss of generality, we suppose that $\underline{xy}=v_8$ and $|N_G(x)\cap \{v_1,v_2,v_6,v_7\}|\geq |N_G(y)\cap \{v_1,v_2,v_6,v_7\}|$. Let $G'=G-\{v_1,v_2,\ldots,v_7,x\}$. Clearly, $G'$ has an $(h_{G'},M)$-coloring $\psi$ by the definition of $G$. Let $\psi'$ be a limitation of $\psi$ on $G'$.
As $x$ has at most two neighbors in $G'$, it follows that
$|L'(x)|\geq M-18M/31=13M/31$. Notice that $|\psi'(y)| \leq 9M/31$. Thus $\psi'$ can be extended to $x$ satisfying $|\psi'(x)\backslash \psi'(y)|\geq 4M/31$.

For $u\in \{v_1,v_2,v_6,v_7\}$, we get $|L'(u)|\geq M-9M/31=22M/31$ as  $u$ has exactly one neighbor in $\{x,y\}$.
Moreover, we may assume that either $|L'(v_2)\setminus L'(v_6)|\geq 4M/31$ or $|L'(v_2)\setminus L'(v_7)|\geq 4M/31$ because of $|\psi'(x)\backslash \psi'(y)|\geq 4M/31$.

 If $|L'(v_2)\setminus L'(v_6)|\geq 4M/31$, then we choose $A\subseteq L'(v_2)\setminus L'(v_6)$ with $|A|=4M/31$ and  $B\subseteq L'(v_1)\setminus A$ with $|B|=5M/31$.
 To define color sets for remaining vertices, we first set  $A \subset \psi'(v_2)$,  $A \subset \psi'(v_5)$,  $B \subset \psi'(v_1)$, and $B \subset \psi'(v_4)$. Then we can complete the extension of $\psi'$ by coloring $v_1,v_2,v_7,v_6,v_5,v_4,v_3$ one by one greedily.
 When we color $v_6$,  note that $|L'(v_6)\setminus (\psi'(v_7)\cup A\cup B)| \geq M-9M/31-8M/31-5M/31=9M/31$ as $A\cap L'(v_6)=\emptyset$.
 When we color $v_3$, notice that $|\cup_{u \in N_G(v_3)} \psi'(u)| \leq 4\cdot (8M/31)-9M/31=23M/31$. Therefore, we can extend $\psi'$  as an $(h_G,M)$-coloring of $G$. This is a contradiction.

If  $|L'(v_2)\setminus L'(v_7)|\geq 4M/31$, then let  $A\subseteq L'(v_2)\setminus L'(v_7)$ with $|A|=4M/31$ and $B\subseteq L'(v_6)\setminus A$ with $|B|=5M/31$. As we did in the previous case, we first set $A \subset \psi'(v_2)$, $A \subset \psi'(v_5)$,
$B \subset \psi'(v_3)$, and $B \subset \psi'(v_6)$. We are also able to color vertices $v_1,v_2,v_3,v_7,v_6,v_5, v_4$ one by one greedily. When we color $v_7$,  note that $|L'(v_7)\setminus (\psi(v_1)\cup A\cup B)| \geq M-9M/31-8M/31-5M/31=9M/31$ as $A\cap L'(v_7)=\emptyset$. Thus $\psi'$ can be extended as an $(h_G,M)$-coloring of $G$, which is a contradiction. \hfill\rule{1ex}{1ex}

\begin{claim}\label{nop52}
$G$ is $P_5^2$-free.
\begin{figure}[htb]
\centering
{\includegraphics[height=0.24\textwidth]{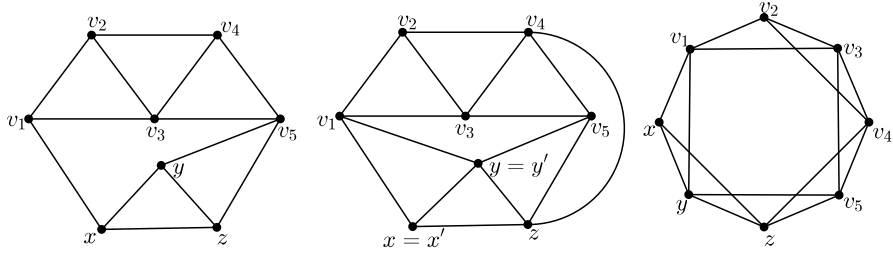}}

\caption{The related graphs in the proof of Claim \ref{nop52}}\label{fig:F11}
\end{figure}
\end{claim}

\noindent
{\bf Proof:} Suppose that $G$ contains a $P_5^2$.  Let $v_1,v_2,\ldots,v_5$ be the vertices of this $P_5^2$.  Since $G$ contains no $K_4$, we have $v_1v_4,v_2v_5\notin E(G)$.

\noindent
{\bf Case 1:} $v_1v_5\in E(G)$.  Let $G'=G-\{v_1,v_2,v_3,v_4,v_5\}$.  Clearly, $G'$ has an  $(h_{G'},M)$-coloring $\psi$ by the definition of $G$. Let $\psi'$ be a limitation of $\psi$ on $G'$.
Note that
$|L'(v_i)|\geq M-9M/31=22M/31$ for $i\in\{1,2,4,5\}$. It follows that $|L'(v_1)\cap L'(v_4)| \geq 13M/31$ and $|L'(v_2)\cap L'(v_5)| \geq 13M/31$. Let  $A\subseteq L'(v_1)\cap L'(v_4)$ with $|A|=6M/31$ and $B\subseteq (L'(v_2)\cap L'(v_5))\setminus A$ with $|B|=7M/31$. We  set $A \subset \psi'(v_1)$, $A \subset \psi'(v_4)$,
$B \subset \psi'(v_2)$ and $B \subset \psi'(v_5)$, which enable us   to color vertices $v_1,v_5,v_2,v_4$ one by one greedily. When we color $v_3$, we note that $|\cup_{u \in N_G(v_3)} \psi(u)| \leq 4\cdot (9M/31) -13M/31=23M/31$. Therefore, we can define $\psi'(v_3)$ properly and  extend $\psi'$  as an $(h_G,M)$-coloring of $G$. This is a contradiction.

\noindent
{\bf Case 2:}  $v_1v_5\notin E(G)$. In this case, we first prove that $d_G(v_i)=4$ for each $i \in \{1,2,4,5\}$ step by step.
 We begin to show $d_G(v_2)=d_G(v_4)=4$.  If  $d_G(v_2)=3$, then let $G'=G-\{v_1,v_2,v_3,v_4,v_5\}$. As the definition of $G$, the graph  $G'$ has an  $(h_{G'},M)$-coloring $\psi$. Let $\psi'$ be a limitation of $\psi$ on $G'$.
Note that $|L'(v_1)|\geq M-9M/31\times 2=13M/31$, $|L'(v_5)|\geq M-9M/31\times 2=13M/31$, and $|L'(v_4)|\geq M-9M/31=22M/31$. It follows that $|L'(v_1)\cap L'(v_4)|\geq 4M/31$. Let  $A\subseteq L'(v_1)\cap L'(v_4)$ with $|A|=4M/31$. We first set $A \subset \psi'(v_1)$ and $A \subset \psi'(v_4)$. Now, we are able to color vertices $v_5,v_1,v_4,v_3$ one by one greedily. When we color $v_2$, we note that $|\cup_{u \in N_G(v_2)} \psi'(u)| \leq 2\cdot (9M/31)+ 8M/31-4M/31=22M/31$. Therefore, we can define $\psi'(v_2)$ properly and  extend $\psi'$  as an $(h_G,M)$-coloring of $G$. This is a contradiction.
Therefore $d_G(v_2)=4$. By the symmetry, we can show $d_G(v_4)=4$.

We next show $d_G(v_1)=d_G(v_5)=4$. If $d_G(v_1)=3$, then let $u_1$ be another neighbor of $v_1$ and $u_2$ be another neighbor of $v_2$. Let $G'=G-\{v_1,v_2,v_3,v_4,v_5\}$.  Then $G'$ has an  $(h_{G'},M)$-coloring $\psi$ by the definition of $G$. Let $\psi'$ be a limitation of $\psi$ on $G'$. If $|\psi'(u_1)\cap \psi'(u_2)|\geq 3M/31$, then we define $A\subseteq \psi'(u_1)\cap \psi'(u_2)$ with $|A|=3M/31$. We  set $A \subset \psi'(v_3)$, which allows us   to color vertices $v_5,v_4,v_3,v_2,v_1$ one by one greedily. This is a contradiction.
If $|\psi'(u_1)\cap \psi'(u_2)|< 3M/31$, then there is a set $A\subseteq \psi'(u_1)\setminus \psi'(u_2)$ with $|A|=3M/31$ and a set $B\subseteq \psi'(u_2)\setminus A$ with $|B|=2M/31$. We first set $A \subset \psi'(v_2)$ and $B \subset \psi'(v_3)$, then we are able to color vertices $v_5,v_4,v_3,v_2,v_1$ one by one greedily. This is a contradiction. Therefore, it follows  that $d_G(v_1)=4$. By the symmetry of $v_1$ and $v_5$, we have $d_G(v_5)=4$. Let us proceed with the assumption $d_G(v_i)=4$ for each $i \in \{1,2,4,5\}$.

Let $G'$ be the graph obtained from $G$ by deleting $v_2,v_3$, and $v_5$ followed by identifying $v_1$ and $v_4$.
Similarly, let $G''$  be the graph obtained from $G$ by deleting $\{v_2,v_3,v_4\}$ followed by identifying $\{v_1,v_5\}$. Claim \ref{cl:c82-v} implies that both $G'$ and $G''$ are $C_8^2$-free. We further assert that one of $G'$ and $G''$ is $K_4$-free. If it is not the case, then the existence of a  $K_4$ in $G''$ implies that there is a triangle $xyz$ outside $P_5^2$ such that $v_1x,v_5y,v_5z \in E(G)$, see the first picture in Figure \ref{fig:F11}. Similarly, a $K_4$ in $G'$ indicates that there exists  a triangle $x'y'z'$ outside $P_5^2$ such that $v_1x',v_1y',v_4z' \in E(G)$. Since  $\Delta(G) \leq 4$,
 we may assume that  $x'=x$. Note that $\{y,z,v_1\}\subseteq N_G(x)$ and $\{y',z',v_1\}\subseteq N_G(x')=N_G(x)$. Thus we have either $y'=y$ or  $y'=z$, say  $y'=y$. In this case, $z'=z$ and $v_4z\in E(G)$ (see the second picture in Figure \ref{fig:F11}), but now $G+v_2x$ contains a $C_8^2$ (see the last picture in Figure \ref{fig:F11}). This is a contradiction to Claim \ref{cl:c82-e} and the assertion follows.

If $G'$ is $K_4$-free, then let $\psi$ be an $(h_{G'},M)$-coloring of $G'$. We may assume that $\psi(v_1)=\psi(v_4)=[1,9M/31]$.
 We will color neighbors of $v_3$ carefully such that each color from $[1,9M/31]$ appears twice in the neighbors of $v_3$. This will ensure that we can color $v_3$ properly.
By Fact 1, there is a limitation $\psi'$ of $\psi$ on $G'$ such that $[1,M/31]\subseteq L'(v_2)$ and $[1,M/31]\subseteq L'(v_5)$. If we define  $A=L'(v_5) \setminus [1,9M/31]$, then it is true that $|A| \geq 4M/31$  as $|L'(v_5)| \geq 13M/31$.  To define $\psi'(v_5)$, we first put
$A \cup [1,M/31] \subset \psi'(v_5)$. In order to make $|\psi'(v_5)|=8M/31$, we  choose at most $3M/31$ colors from $[1+M/31,9M/31]$ and add them to $\psi'(v_5)$.
To color $v_4$, we set $([1+M/31,9M/31]\setminus \psi'(v_5))\subseteq \psi'(v_4)$ and add extra colors greedily
such that $[1,M/31]\cap \psi'(v_4)=\emptyset$.
 Define $\psi'(v_1)=[1+M/31,9M/31]$. Note that $|\psi'(v_1)\cap \psi'(v_4)| \geq 3M/31$. It follows that $\psi'$ can be extended to $v_2$ such that $[1,M/31]\subseteq \psi'(v_2)$.  Since each color in $[1,9M/31]$ appears twice in neighbors of $v_3$,  then $|L'(v_3)|\geq  M-(4\times8M/31-9M/31)=8M/31$. Therefore, we can define $\psi'(v_3)$ properly to get an  $(h_G,M)$-coloring of $G$, a contradiction.

 If $G''$ is $K_4$-free, then we next show that $G$ has an $(h_G,M)$-coloring.   Let $\psi$ be an $(h_{G''},M)$-coloring of $G''$. We may assume that $[1,8M/31] \subseteq \psi(v_1)=\psi(v_5)$. Consequently, $[1,8M/31] \subseteq L(v_1)$ and $[1,8M/31] \subseteq L(v_5)$. Note that $|L(v_1)| \geq 11M/31$, which enables us to choose $A\subseteq L(v_1)$ such that $|A|=M/31$ and $A \cap [1,8M/31]=\emptyset$. By Fact 1, we can define $\psi'$ as a limitation of $G-\{v_1,v_2,v_3,v_4,v_5\}$ on $\psi$ such that $A \subseteq L'(v_4)$.
   Observe that $|L'(v_2)|\geq M-9M/31=22M/31$ and $|L'(v_5)|\geq M-2\times 9M/31=13M/31$. Thus we get  $|L'(v_2)\cap L'(v_5)|\geq 4M/31$. Let $B\subseteq L'(v_2)\cap L'(v_5)$ such that  $|B|=2M/31$ and $A\cap B=\emptyset$.
 To color $v_5$, we first put $B$ in  $\psi'(v_5)$ and add extra colors from $[1,8M/31] \setminus B$ to $\psi'(v_5)$ so that
 $|\psi'(v_5)|=8M/31$. To color $v_1$, we put $(\psi'(v_5)\setminus B)\cup A$ in $\psi'(v_1)$ followed by adding extra colors to  $\psi'(v_1)$ so that  $|\psi'(v_1)|=8M/31$ and
 $\psi'(v_1)\cap B=\emptyset$. Extend $\psi'$ to $v_2$ such that $B\subseteq \psi'(v_2)$ and $ \psi'(v_2) \cap A=\emptyset$.
Since $B\subseteq \psi'(v_2)\cap \psi'(v_5)$, we can  define $\psi'(v_4)$ properly with  $A\subseteq \psi'(v_4)$.
  Since each color in $A \cup \psi'(v_5)$ appears twice in neighbors of $v_3$,  then $|L'(v_3)|\geq  M-(4\times8M/31-9M/31)=8M/31$. We can extend  $\psi'$ to $v_3$ to get an  $(h_G,M)$-coloring of $G$, a contradiction.
   \hfill\rule{1ex}{1ex}

~~

Let $F_4$ be a graph consisting of a triangle $uvw$ and a vertex $x$ adjacent to $u$ and $v$. We will keep the labelling of vertices of this graph in the rest of this paper.

\begin{claim}\label{cl:db-triangle-3-3-1}
$G$ contains no $F_4$ such that $d_G(u)=d_G(v)=d_G(w)=3$.
\end{claim}

\noindent
{Proof:} 
Let $G'=G-\{u,v,w,x\}$. Then $G'$ has an $(h_{G'},M)$-coloring $\psi$. Let $\psi'$ be a  limitation of $\psi$ on $G'$.
As  $\delta(G) \geq 3$, it follows that either $|\psi'(y)|=8M/31$ or $|\psi'(y)|=9M/31$ for each $y \in V(G')$. Note that $|L'(w)| \geq 22M/31$ and $|L'(x)| \geq M-(d_G(x)-2)\times9M/31=(49-9d_G(x))M/31$.
 Then $|L'(w) \cap L'(x)|\geq 22M/31+(49-9d_G(x))M/31-M=(40-9d_G(x))M/31\geq (8-d_G(x))M/31$.
 We can first color $w$ and $x$ such that $|\psi'(x) \cap \psi'(w)|=(8-d_G(x))M/31$. As $|\psi'(x) \cup \psi'(w)|=9M/31+(12-d_G(x))M/31-(8-d_G(x))M/31=13M/31$, we can define the coloring of $u$ and $v$ greedily. Thus $\psi'$ can be extended as an $(h_G,M)$-coloring of $G$, which is a contradiction and the claim follows.
  \hfill\rule{1ex}{1ex}

\begin{claim} \label{noG1}
 $G$ does not contain $F_7$ with $d_G(u)=3$ as a subgraph, where     $F_7$ is shown in Figure \ref{fig:F7}.
\begin{figure}[!htb]
\centering
{\includegraphics[height=0.25\textwidth]{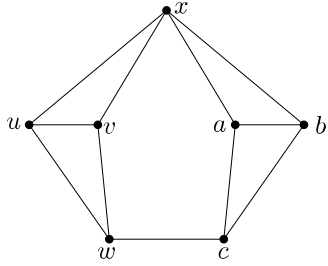}}
\caption{The graph $F_7$}\label{fig:F7}
\end{figure}

\end{claim}
\noindent
{\bf Proof:} Suppose that $G$ contains $F_7$ as a subgraph. By Claim \ref{nop52}, $F_7$ is an induced subgraph of $G$, and by Claim \ref{cl:db-triangle-3-3-1},  $d_G(a)+d_G(b)+d_G(c)\geq 3+3+4=10$. We claim that there is a set coloring $\psi'$ of $G-\{u,v\}$ such that $|\psi'(y)|=h_G(y) M$ for each $y \in V(G)-\{u,v\}$ and $|\psi'(x) \cap \psi'(w)| \geq 3M/31$.

Let $G'=G-\{u,v,w,x\}$. As $G$ is a minimum counterexample, $G'$ has an $(h_{G'},M)$-coloring $\psi$.
For each $y \in \{a,b,c\}$, the fact $d_{G'}(y) \leq d_G(y)-1$ implies that $h_{G'}(y) \geq h_G(y) +1/31$.
There are two cases depending on the degree of $w$.

\noindent
{\bf Case 1:} $d_G(w)=4$.  As $\{u,v,c\} \subset N_G(w)$, we assume that $w'$ is another neighbor of $w$.
The assumption $d_G(a)+d_G(b)+d_G(c)\geq 10$ implies  that $d_{G'}(a)+d_{G'}(b)+d_{G'}(c)\geq 7$ and $|\psi(a) \cup \psi(b) \cup \psi(c)|\leq 29M/31$. Let $A\subseteq M \setminus (\psi(a) \cup \psi(b) \cup \psi(c))$ with $|A|=2M/31$. Denote $B=A\cap\psi(w')$ and $C= A\setminus B$. Note that $|(\psi(a)\cup \psi(b))\setminus\psi(w')|\geq 8M/31$.

If $|B|\leq M/31$, then let $D\subseteq (\psi(a) \cup \psi(b)) \setminus\psi(w')$ such that $|D|=M/31$.
By Fact \ref{fact1}, we can define a limitation $\psi'$ of $\psi$ on $G'$ such that $B\subseteq L'(w)$ and $D\subseteq L'(x)$.
Set $(A\cup D) \subset \psi'(w)$ and $(A\cup D) \subset \psi'(x)$. Then we pick extra colors greedily for $\psi'(w)$ and $\psi'(x)$ to satisfy
$|\psi'(w)|=Mh_G(w)$ and $|\psi'(x)|=Mh_G(x)$. We note that $|\psi'(x) \cap \psi'(w)| \geq |A|+|D|=3M/31$.

If $|B|> M/31$, then let $B'\subseteq B$ such that $|B'|= M/31$. If there is a vertex $y \in \{a,b\}$ such that $|\psi(y) \setminus \psi(w')|<M/31$, say $y=a$ by the symmetry, then there are sets $C \subset \psi(b) \setminus \psi(w')$ and $D \subset \psi(c) \setminus \psi(w')$ with $|C|=|D|=M/31$ as $\psi(a)$, $\psi(b)$, and $\psi(c)$ are pairwise disjoint.
To define a limitation $\psi'$ of $\psi$ on $G'$, we remove colors of $C$ from $\psi(b)$, colors of $D$ from $\psi(c)$, and colors of $B'$ from $\psi(w')$. Thus we have $B' \cup C \cup D  \subset L'(w)$ and  $B' \cup C \cup D \subset L'(x)$.
 Set $B' \cup C \cup D \subset \psi'(w)$ and $B' \cup C \cup D \subset \psi'(x)$ in this case. If
$|\psi(a) \setminus \psi(w')| \geq M/31$ and $|\psi(b) \setminus \psi(w')| \geq M/31$, then let $C \subset \psi(a) \setminus \psi(w')$ and $D \subset \psi(b) \setminus \psi(w')$ with $|C|=|D|=M/31$.   We can define a limitation $\psi'$ of $\psi$ on $G'$ such that $B'   \subset L'(w)$ by Fact \ref{fact1}. Moreover, we can assume that  $C \cup D \subset L'(x)$ as we can remove colors of $C$ from $\psi(a)$ and colors of $D$ from $\psi(b)$ when we define the limitation. Similarly, set  $B' \cup C \cup D \subset \psi'(w)$ and $B' \cup C \cup D \subset \psi'(x)$. Then we pick extra colors greedily for $\psi'(w)$ and $\psi'(x)$ to satisfy
$|\psi'(w)|=Mh_G(w)$ and $|\psi'(x)|=Mh_G(x)$. We note that $|\psi'(x) \cap \psi'(w)| \geq |B'|+|C|+|D|=3M/31$.

\noindent
{\bf Case 2:} $d_G(w)=3$. Let $\psi'$ be a limitation of $\psi$ on $G-\{u,v,w,x\}$.  As $|\psi'(a) \cup \psi'(b) \cup \psi'(c)|\leq 27M/31$, we are able to choose $A \subseteq M \setminus (\psi'(a) \cup \psi'(b) \cup \psi'(c))$ with $|A|=3M/31$. Set $A \subseteq \psi'(w)$ and $A \subseteq \psi'(x)$. Then we pick extra colors greedily for $\psi'(w)$ and $\psi'(x)$ to satisfy
$|\psi'(w)|=Mh_G(w)$ and $|\psi'(x)|=Mh_G(x)$. We note that $|\psi'(x) \cap \psi'(w)| \geq 3M/31$.

To finish the coloring,  Claim \ref{cl:db-triangle-3-3-1} implies that $d_G(w)+d_G(v) \geq 7$. We first color $v$ greedily. Note that $|\psi'(x) \cup \psi'(v) \cup \psi'(w)| \leq 2\cdot 8M/31+9M/31-3M/31=22M/31$ in each case. Thus we can color $u$ properly and get an $(h_G,M)$-coloring of $G$, which is a contradiction.  \hfill\rule{1ex}{1ex}

\begin{claim}
\label{cl:db-triangle-3-3}
 $G$ contains no $F_4$ such that $d_G(u)=3$.
\end{claim}

\noindent{\bf Proof:}  Suppose that $G$ contains an $F_4$ such that $d_G(u)=3$.
Let $G'$ be the graph obtained from
$G-\{u,v\}$ by identifying $w$ and $x$. Claim \ref{cl:c82-v} tells us that $G'$ contains no $C_8^2$. If $G'$ contains a $K_4$, then $G$ contains an $F_7$, which is a contradiction to Claim \ref{noG1}. Therefore, $G'$ is $K_4$-free and $G'$ has an $(h_{G'},M)$-coloring $\psi$. Note that $\psi(w)=\psi(x)$.
We can extend $\psi$ as an $(h_G,M)$-coloring of $G$ by coloring $v$ and $u$ greedily, which is a contradiction.  \hfill\rule{1ex}{1ex}

~~~

Let $F_5$ be a graph consisting of an edge $uv$ and three common neighbors $a$, $b$, and $c$. 

\begin{claim}
\label{cl:tr-triangle}
 $G$ contains no $F_5$.
\end{claim}
\noindent{\bf Proof:} Suppose  that $G$ contains an $F_5$. As $G$ is $K_4$-free, we get that $\{a,b,c\}$ is an independent set.
 There are three cases depending on the number of three vertices from $\{a,b,c\}$.


\noindent
{\bf Case 1:} There are two vertices from $\{a,b,c\}$ having degree three. We can  assume that $d_G(a)=d_G(b)=3$ and $3\leq d_G(c)\leq 4$.  Let $G'=G-\{u,v,a,b,c\}$. Then  $G'$ has an $(h_{G'},M)$-coloring $\psi$. Let $\psi'$ be a limitation of $\psi$ on
 $G-\{u,v,a,b,c\}$.
Note that $|L'(a)|\geq 22M/31$, $|L'(b)|\geq 22M/31$, and $|L'(c)|\geq 13M/31$. Thus $|L'(a)\cap L'(b)|\geq 13M/31$, $|L'(a)\cap L'(c)| \geq 4M/31$, and $|L'(b)\cap L'(c)| \geq 4M/31$.
Let $A\subseteq L'(a)\cap L'(c)$ with $|A|=4M/31$ and  $B\subseteq L'(b)\cap L'(c)$ such that $|B|=4M/31$. Moreover,
we choose $C\subseteq L(a)\cap L(b)$ such that $|C|=5M/31$ and $C\cap(A\cup B)=\emptyset$. Set $A\cup C \subset \psi'(a)$, $B \cup C \subset \psi'(b)$, and $A \cup B \subset \psi'(c)$.
Then $\psi'$ can be extended as an $(h_G,M)$-coloring of  $G$ by choosing colors for $a,b,c,u,v$ one by one greedily, which is a contradiction.

\noindent
{\bf Case 2:}   
There is exactly one vertex of $\{a,b,c\}$ with degree three. Then we assume that $d_G(a)=3$ and $d_G(b)=d_G(c)=4$.
Let $G'$ (and $G''$) be the graph obtained from $G-\{u,v,c\}$ (and $\{u,v,b\}$) by identifying $a$ and $b$ (by identifying $a$ and $c$) respectively. Note that both $G'$ and $G''$ are $C_8^2$-free by Claim \ref{cl:c82-v}.  There are two subcases.

\vspace{0.2cm}
\noindent
{\bf Subcase 2.1:} One of $G'$ and $G''$ is $K_4$-free. By the symmetry, we can assume that $G'$ is $K_4$-free. Thus $G'$ has an $(h_{G'},M)$-coloring $\psi$. Note that $\underline{ab}$ is a 3-vertex of $G'$. Then we may assume that  $\psi(\underline{ab})=[1,9M/31]$.  Let $A=[1,M/31]$. By Fact 1, we can let $\psi'$ be a limitation of $\psi$ on $G-\{u,v,a,b,c\}$ such that $A\subseteq L'(c)$. Note that $[1,9M/31] \subseteq L'(a) \cap L'(b)$, $|L'(a)| \geq 22M/31$, and $|L'(c)| \geq 13M/31$.  It follows that $|L'(a) \cap L'(c)|\geq 4M/31$. Choose $B \subset L'(a) \cap L'(c)$ such that $|B|=M/31$ and $A\cap B=\emptyset$.  Define $\psi'(b)=[1,8M/31]$. We put $[1,8M/31]\cup B\subseteq \psi'(a)$ and $A\cup B \subseteq \psi'(c)$ followed by adding extra colors to $\psi'(a)$  and $\psi'(c)$ greedily.  We can construct an $(h_G,M)$-coloring of $G$ by choosing colors for $u$ and $v$ properly. A contradiction.

\begin{figure}[!htb]
\centering
{\includegraphics[height=0.3\textwidth]{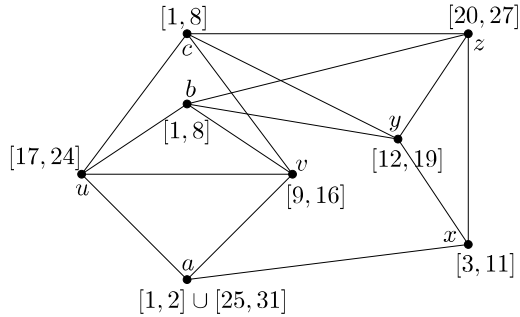}}

\caption{The graph $F_8$ and an $(h_{F_8}, 31)$-coloring of $F_8$}\label{fig:F8}
\end{figure}

\vspace{0.2cm}
\noindent
{\bf Subcase 2.2:} Both $G'$ and $G''$ contain a $K_4$.  In this case, by the symmetry,  there are vertices $x,y,$ and $z$ outside $F_5$ such that $\{xy,xz,yz,ax,yb,yc,zb,zc\} \subset E(G)$.
Since $\Delta(G)\leq 4$ and $G$ is 2-connected, we have $G$ is the graph $F_8$ as shown in Figure \ref{fig:F8}. It is easy to see that $G$ has an $(h_{G},31)$-coloring. A contradiction.

\noindent
{\bf Case 3:}  $d_G(a)=d_G(b)=d_G(c)=4$.  Let $G'$ be the graph obtained from $G$ by deleting vertices from $\{u,v,c\}$ followed by identifying vertices from $\{a,b\}$. Similarly, let $G''$ (and $G'''$) be the graph obtained from $G$ by deleting vertices from $\{u,v,b\}$ (and $\{u,v,a\}$) followed by identifying vertices from $\{a,c\}$ (and vertices from $\{b,c\}$) respectively. Note that Claim \ref{cl:c82-v} implies that all of them are $C_8^2$-free.
 We proceed to complete the proof with the assumption that
 one of $G'$, $G''$, and $G'''$ is $K_4$-free.  By the symmetry, we suppose that $G'$ is $K_4$-free. Let $\psi$ be an $(h_{G'},M)$-coloring of $G'$. We  assume that $\psi'(a)=\psi'(b)=[1,8M/31]$ and fix $C \subset [1,8M/31]$ with $|C|=M/31$.
 By Fact 1, we can
let $\psi'$ be a limitation of $\psi$ on $G-\{u,v,a,b,c\}$ such that $C\subseteq L'(c)$.
Set $C \subset \psi'(c)$ and add extra colors greedily. Now, $|\psi'(a) \cup \psi'(b) \cup \psi'(c)| \leq 8M/31+7M/31=15M/31$ and we can choose colors for $u$ and $v$ greedily to get an $(h_G,M)$-coloring of $G$, a contradiction.

It remains to show the assertion that one of $G'$, $G''$, and $G'''$ is $K_4$-free. If $G'$ contains a $K_4$, then there is a triangle with vertex set $\{x,y,z\}$ outside $F_5$. By the symmetry, we can assume that $\{ax,by,bz\} \subset E(G)$.
Since $G[\{a,b,x,y,z\}]$ is not a $P_5^2$ by Claim \ref{nop52}, we have $ay,az\notin E(G)$.
If $cx\in E(G)$, then $G''$ is $K_4$-free as $G$ is $P_5^2$-free. Now we assume that $cx\notin E(G)$.
The existence of a $K_4$ in $G'''$  implies that there exists $x'$ such that $cx',x'y,x'z\in E(G)$. Now we can observe that $G''$ is $K_4$-free. The assertion is proved and the claim follows.
\hfill\rule{1ex}{1ex}

\begin{claim}
\label{cl:triangle-3-3-3}
 $G$ contains no triangle $uvw$ such that $d_G(u)=d_G(v)=d_G(w)=3$.
\end{claim}

\noindent{\bf Proof:}   Suppose that $G$ contains a triangle $uvw$ with $d_G(u)=d_G(v)=d_G(w)=3$. Let $u',v',$ and $w'$ be the  neighbors of $u,v,w$ which is  not on the triangle, respectively. Claim \ref{cl:db-triangle-3-3} implies that $u',v'$, and $w'$ are distinct. Let $G'=G-\{u,v,w\}$.
  By Claim \ref{cl:c82-e}, $G'+u'v'$, $G'+u'w'$, and  $G'+v'w'$ are $C_8^2$-free.
We claim that either $G'+u'v'$, or $G'+u'w'$, or $G'+v'w'$ is $K_4$-free. Otherwise, either $G$ contains a $K_4$, or there exists $xy\in E(G')$ such that $\{x,y\}$ together with each of  $u',v'$, and $w'$ form a triangle, which is a contradiction to Claim \ref{cl:tr-triangle}.

Assume  by the symmetry that $G'+u'v'$ is $K_4$-free. As the definition of $G$, the graph $G'+u'v'$ has an $(h_{G'},M)$-coloring $\psi$. Since $\psi(u')\cap\psi(v')=\emptyset$, any limitation $\psi'$ of $\psi$ on $G'$ can be extended as an $(h_G,M)$-coloring of $G$ by Fact \ref{fact2}, a contradiction.
 \hfill\rule{1ex}{1ex}
%
%

\begin{claim}
\label{cl:path-3-3-3}
 $G$ contains no path $u_1u_2u_3$ such that $d_G(u_1)=d_G(u_2)=d_G(u_3)=3$.
\end{claim}

\noindent{\bf Proof:} Suppose  that $G$ contains a path $u_1u_2u_3$ such that $d_G(u_1)=d_G(u_2)=d_G(u_3)=3$. By Claim \ref{cl:triangle-3-3-3}, we have $u_1u_3\notin E(G)$.
For $i\in \{1,3\}$, we assume that $u_i'$ and $u_i''$ are the other two neighbors of $u_i$.

If $d_G(u_1')+d_G(u_1'')=d_G(u_3')+d_G(u_3'')=6$, then let $G'$ be the graph obtained from $G$ by deleting $u_2$ followed by identifying $u_1$ and $u_3$. As $d_G(u_1')+d_G(u_1'')=d_G(u_3')+d_G(u_3'')=6$ and Claim \ref{cl:triangle-3-3-3}, the graph $G'$ is $K_4$-free. Claim \ref{cl:c82-v} implies that $G'$ is also $C_8^2$-free. Thus $G'$ has an $(h_{G'},M)$-coloring $\psi$. We can easily extend $\psi$ as an $(h_G,M)$-coloring of $G$, a contradiction.

If either $d_G(u_1')+d_G(u_1'') \geq 7$ or $d_G(u_3')+d_G(u_3'') \geq 7$, then by the symmetry, we assume that $d_G(u_1')+d_G(u_1'') \geq 7$. Let $G'=G-\{u_1,u_2,u_3\}$. Then $G'$ has an $(h_{G'},M)$-coloring $\psi$ and let $\psi'$ be a limitation of $\psi$ on $G-\{u_1,u_2,u_3\}$.
Note that $|L'(u_1)| \geq 14M/31$, $|L'(u_2)| \geq 22M/31$, and $|L'(u_3)|\geq 13M/31$, here we applied the assumption $d_G(u_1')+d_G(u_1'') \geq 7$. Let $A\subseteq L'(u_2)\setminus L'(u_3)$ such that
$|A|=5M/31$. We first extend $\psi'$ to $u_1$ such that $\psi'(u_1)\cap A=\emptyset$ followed by extending $\psi'$ to $u_2$ such that $A\subseteq \psi'(u_2)$. Since $|L'(u_3)\setminus \psi'(u_2)|\geq 9M/31$,  we can extend $\psi'$  to $u_3$ and get an $(h_G,M)$-coloring of $G$, a contradiction.
   \hfill\rule{1ex}{1ex}

\begin{claim}
\label{cl:triangle-3-3}
 $G$ contains no triangle $uvw$ such that $d_G(u)=d_G(v)=3$.
\end{claim}

\noindent{\bf Proof:} Suppose  that $G$ contains a triangle $uvw$ such that $d_G(u)=d_G(v)=3$. Then $d_G(w)=4$ by Claim \ref{cl:triangle-3-3-3}. Let $u'$ and $v'$ be the neighbors of $u$ and $v$ not on the triangle, respectively.  Claim \ref{cl:db-triangle-3-3} gives us that $u'\neq v'$, and $wu',wv'\notin E(G)$. Claim \ref{cl:path-3-3-3} implies that $d_G(u')=d_G(v')=4$.

If $u'v' \in E(G)$, then let $G'=G-\{u,v\}$. Clearly, $G'$ has an $(h_{G'},M)$-coloring $\psi$, where $\psi(u')\cap\psi(v')=\emptyset$. Thus a limitation of $\psi$ on $G'$ can be extended as an $(h_G,M)$-coloring of $G$ by Fact \ref{fact2}, a contradiction. We need only to consider the case where $u'v' \notin E(G)$.
Let $G'=G-\{u,v\}+u'v'$. By Claim \ref{cl:c82-e}, $G'$ is $C_8^2$-free. If $G'$ is also $K_4$-free, then $G'$ has an $(h_{G'},M)$-coloring $\psi$. Since $\psi(u')\cap\psi(v')=\emptyset$, a limitation of $\psi$ on $G'$ can be extended as an $(h_G,M)$-coloring of $G$ by Fact \ref{fact2}, a contradiction.
 Therefore, $G'$ contains a $K_4$, which implies that  there exist $xy\in E(G)$ such that $u'x,u'y,v'x,v'y\in E(G)$. If $u'w\in E(G')$, then $G$ contains $F_7$ as a subgraph, which is a contradiction to Claim \ref{noG1}. Thus $u'w\notin E(G')$. By the symmetry, we also have $v'w \notin E(G')$.

We first assume that $wx\in E(G)$.  Claim \ref{cl:tr-triangle} gives that $wy\notin E(G)$. Let $w'$ be the other neighbor of $w$. Without raising any confusion, we define  $G'=G-\{u,v,w\}+v'w'$. By Claim \ref{cl:c82-e}, $G'$ is $C_8^2$-free. Note that  $G'$ is also  $K_4$-free by checking degrees of neighbors of $v'$. Thus $G'$ has an $(h_{G'},M)$-coloring $\psi$. Let $\psi'$ be a limitation of $\psi$ on $G-\{u,v,w\}$.
If $|\psi'(u')\setminus \psi'(v')|\geq 4M/31$, then $\psi'$ can be extended as an $(h_G,M)$-coloring of $G$ by Fact \ref{fact2}, a contradiction.
  Therefore $|\psi'(u') \setminus \psi'(v')| < 4M/31$, which implies $|\psi'(u') \cap \psi'(v')|>4M/31$.  Let $B=\psi'(u') \cap \psi'(v')$ with $|B|=4M/31$. We first put $B \subset \psi'(w)$ and add extra colors greedily to ensure $|\psi'(w)|=8M/31$. This is possible because $u'x \in E(G)$ and $w'v' \in E(G')$.   Then we can define $\psi'(u)$ and $\psi'(v)$ one by one greedily. In this case, we can construct an $(h_G,M)$-coloring of $G$, which is a contradiction.

Now we can  assume that $wx\notin E(G)$. By the symmetry, we have $wy\notin E(G)$.
Let $w_1$ and $w_2$ be the neighbors of $w$ not on the triangle.
Define  $G'=G-\{u,v,w\}+w_1u'+w_2v'$.
We first claim that $G'$ is $K_4$-free. Otherwise, this $K_4$ must contain exactly one of the edges $\{w_1u',w_2v'\}$ as $u'v'\notin E(G)$. We can assume  that $K_4$ contains the edge $w_1u'$.
As $\Delta(G) \leq 4$, the only possibility for the vertex set of this $K_4$ is $\{w_1,u',x,y\}$. Now $x$ and $y$ have three common neighbors $u',v'$ and $w_1$, which is a contradiction to Claim \ref{cl:tr-triangle}.

If $G'$ contains no $C_8^2$, then $G'$ has an $(h_{G'},M)$-coloring $\psi$. Let $\psi'$ be a limitation of $\psi$ on $G-\{u,v,w\}$. If $|\psi'(v') \cap \psi'(u')| \geq 3M/31$, then let $A \subset \psi'(v') \cap \psi'(u')$ with $|A|=3M/31$. This is possible as $w_1u',w_2v' \in E(G')$.  We put $A$ in $\psi'(w)$ and add extra colors to $\psi'(w)$ such that $|\psi'(w)|=8M/31$. Then $\psi'$ can be extended as an $(h_G,M)$-coloring of $G$ by defining $\psi'(u)$ and $\psi'(v)$ one by one greedily.  If $|\psi'(v') \cap \psi'(u')| <3M/31$, then $|\psi'(v') \setminus \psi'(u')| > 5M/31$. 
By Fact \ref{fact2}, $\psi'$ can be extended to get an $(h_G,M)$-coloring of $G$, a contradiction.

\begin{figure}[!htb]
\centering
{\includegraphics[height=0.3
\textwidth]{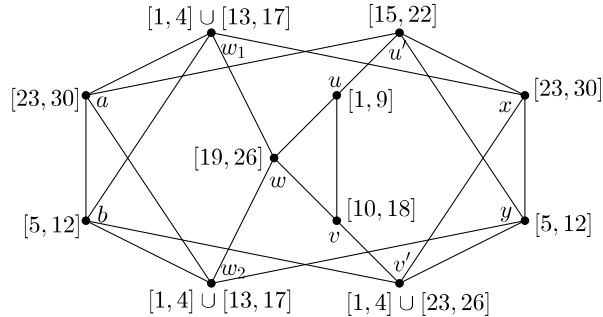}}
\caption{The graph $F_{11}$ and an $(h_{F_{11}},31)$-coloring of $F_{11}$}\label{fig:11}
\end{figure}

Now we assume that $G'$ contains a $C_8^2$.
By Claim \ref{cl:c82-e}, this $C_8^2$ must contain  edges $w_1u'$ and $w_2v'$.
 We assume that $N_{G}(u')=\{a,w_1,x,y\}$ and $N_{G}(v')=\{b,w_2,x,y\}$. The fact that all vertices in $G'$ have degree at most four implies that $G'$ is the whole graph $C_8^2$.
Since $G[\{a,u',v',x,y\}]$ contains no $P_5^2$ as a subgraph by Claim \ref{nop52}, it follows that $ax,ay\notin E(G)$. Then we get  $aw_1,aw_2,ab\in E(G)$. By the symmetry of $a$ and $b$, we have $bw_1,bw_2\in E(G)$. Therefore, $G$ must be the graph $F_{11}$ as shown in Figure \ref{fig:11}. However,  $F_{11}$ has an $(h_{F_{11}},31)$-coloring as shown by Figure \ref{fig:11}, a contradiction.

 The proof of the claim is complete.
\hfill\rule{1ex}{1ex}

\begin{claim}
\label{cl:path 3-3-4-3}
 $G$ contains no  path $u_1u_2u_3u_4$ such that $d_G(u_1)=d_G(u_2)=d_G(u_4)=3$.
\end{claim}

\noindent{\bf Proof:} Suppose that $G$ contains a path $u_1u_2u_3u_4$ such that $d_G(u_1)=d_G(u_2)=d_G(u_4)=3$.  Claim \ref{cl:triangle-3-3} implies that $u_1u_2u_3u_4$ is an induced path.
By Claim \ref{cl:path-3-3-3}, $d_G(u_3)=4$ while other
 neighbors of $u_1$ and $u_2$ not in this path are 4-vertices.

 Let $G'=G-\{u_1,u_2\}$.  Clearly, $G'$ has an $(h_{G'},M)$-coloring $\psi$.
Suppose that $\psi'$ is a  limitation of $\psi$ on $G-\{u_1,u_2,u_3,u_4\}$.
Note that $|L'(u_1)|=15M/31$, $|L'(u_2)|=23M/31$,  $|L'(u_3)|\geq 13M/31$, and $|L'(u_4)|\geq 13M/31$. Moreover, $|L'(u_3)\cup L'(u_4)|\geq 17M/31$
 as it contains $\psi(u_3) \cup \psi(u_4)$ as a subset. By deleting colors not in $\psi(u_3) \cup \psi(u_4)$ if it necessary, we may assume that $|L'(u_3)|=|L'(u_4)|=13M/31$. Then
 $|L'(u_3)\setminus L'(u_4)|\geq 4M/31$ because of  $|L'(u_3)\cup L'(u_4)|\geq 17M/31$. Let $A\subseteq L'(u_3)\setminus L'(u_4)$ and $B\subseteq L'(u_2)\setminus L'(u_1)$ such that   $|A|=4M/31$, $|B|=3M/31$, and $A\cap B=\emptyset$.  Put $A \subset \psi'(u_3)$ and $B \subset \psi'(u_2)$. Choosing extra colors for $u_3$, $u_2$, $u_1$ and $u_4$ one by one greedily, we obtain an $(h_G,M)$-coloring of $G$, a contradiction.
\hfill\rule{1ex}{1ex}

\begin{claim}
\label{cl:db-triangle-w-x}
 $G$ contains no $F_4$ such that $d_G(w)=d_G(x)=3$.
\end{claim}

\noindent{\bf Proof:}  Suppose  that $G$ contains an $F_4$ such that $d_G(w)=d_G(x)=3$.
Claim \ref{cl:triangle-3-3} indicates that $d_G(u)=d_G(v)=4$.
For each $z \in \{u,v,w,x\}$, let $z'$ be the neighbor of $z$ which is not in $F_4$. Claim \ref{cl:path 3-3-4-3} implies that $d_G(w')=d_G(x')=4$. Moreover, $u',v',$ and $w'$ are all distinct by Claim \ref{nop52} and Claim \ref{cl:tr-triangle}. Similarly,
 $u',v',$ and $x'$ are  distinct.
Let $G'=G-\{u,v,w,x\}$. Then $G'$ is $C_8^2$-free and $K_4$-free. Thus $G'$ has an $(h_{G'},M)$-coloring $\psi$. There are two cases depending on the size of $\psi(w')\cap \psi(x')$.

\noindent
{\bf Case 1:} $|\psi(w')\cap \psi(x')|\geq 2M/31$. Let $A\subseteq(\psi(w')\cap \psi(x'))$ and $B\subseteq(\psi(w')\cap \psi(x'))$ such that $A\cap B=\emptyset$ and $|A|=|B|=M/31$.
 By Fact 1, there exists a limitation $\psi'$ of $\psi$ on $G-\{u,v,w,x\}$ such that $A\subseteq L'(u)$ and  $B\subseteq L'(v)$.
  Therefore, it is proper to put $A$ in $\psi'(u)$ and $B$ in $\psi'(v)$. We add extra colors to $\psi'(u)$ and $\psi'(v)$ such that $|\psi'(u)|=|\psi'(v)|=8M/31$. Then $\psi'$ can be extended as an $(h_G,M)$-coloring of $G$ by picking colors for $\psi'(w)$ and $\psi'(x)$ one by one greedily. Thus we get an $(h_G,M)$-coloring of $G$, which is a contradiction.

\noindent
{\bf Case 2:} $|\psi(w')\cap \psi(x')|< 2M/31$. We have $w' \not = x'$ in this case. Note that $|\psi(w')\cup \psi(x')|> 16M/31$ and $|L(u)|\geq 21M/31$. It follows that  $|L(u)\cap (\psi(w')\cup \psi(x'))|>6M/31$. By the symmetry, we may assume that $|L(u) \cap \psi(w')|>3M/31$.
Let $A \subset \psi(x')$ and $B \subset \psi(x')$ such that $A \cap B=\emptyset$ and $|A|=|B|=M/31$. Additionally, we choose $C \subset L(u) \cap \psi(w')$ with $|C|=3M/31$.
Fact \ref{fact1} gives a limitation $\psi'$ of $\psi$ on $G-\{u,v,w,x\}$ with the property $A \subset L'(u)$ and $B \subset L'(v)$. Furthermore, we have $C \subset L'(u) \cap \psi'(w')$ by deleting colors from $\psi(w')$ properly.  To color $u$ and $v$, we put $A \cup C$ in $\psi'(u)$ and $B$ in $\psi'(v)$ followed by adding extra colors greedily. As the selection of sets $A,B,$ and $C$, we are able to color $x$ and $w$ one by one, and we obtained an $(h_G,M)$-coloring of $G$. A contradiction.
\hfill\rule{1ex}{1ex}

\begin{claim}
\label{cl:db-triangle-x}
 $G$ contains no $F_4$ such that $d_G(x)=3$.
\end{claim}

\noindent{\bf Proof:}  Suppose  that $G$ contains an $F_4$ such that $d_G(x)=3$.  Claim \ref{cl:triangle-3-3} and Claim \ref{cl:db-triangle-w-x} yield that $d_G(u)=d_G(v)=d_G(w)=4$.
For each $y \in \{u,v,x\}$, let $y'$ be another neighbor of $y$ not in $F_4$.  By Claim \ref{nop52} and Claim \ref{cl:tr-triangle}, $u',v',$ and $x'$ are distinct.
Similarly, if $w'$ and $w''$ are  neighbors of $w$ not in $F_4$, then $u',v',w',$ and $w''$ are distinct.

 It is easy to see that either $G-\{u,v,w,x\}+u'x'$ or $G-\{u,v,w,x\}+v'x'$ is $K_4$-free.  Assume  that $G'=G-\{u,v,w,x\}+u'x'$ is $K_4$-free.
 Claim \ref{cl:c82-e} implies that $G'$ is also $C_8^2$-free. Thus $G'$  has an $(h_{G'},M)$-coloring $\psi$. Let  $A\subseteq \psi(u')$ and $B\subseteq \psi(u')$ such that $A\cap B=\emptyset$ and $|A|=|B|=M/31$, and let $C\subseteq \psi(x')$ with $|C|=3M/31$. Note that $A,B$, and $C$ are pairwise disjoint.
By Fact \ref{fact1}, we may assume that $\psi'$ is a  limitation of $\psi$ on $G-\{u,v,w,x\}$ such that $A\subseteq L'(w)$ and $B\subseteq L'(v)$. Put $A$ in $\psi'(w)$, $B$ in $\psi'(v)$, and $C$ in $\psi'(u)$. Choosing extra colors for $w$, $v$, $u$ and $x$ one by one greedily, we obtain an $(h_G,M)$-coloring of $G$, a contradiction.
 \hfill\rule{1ex}{1ex}

\begin{claim}
\label{cl:triangle-3}
 $G$ contains no triangle $uvw$ such that $d_G(u)=3$.
\end{claim}

\noindent{\bf Proof:} Suppose  that $G$ contains a triangle $uvw$ such that $d_G(u)=3$. By Claim \ref{cl:triangle-3-3}, we have $d_G(v)=d_G(w)=4$.
Let $x$ be the neighbor of $u$ not on the triangle.  Claim \ref{cl:db-triangle-3-3} implies that $vx,wx\notin E(G)$. Suppose that $N_G(w)=\{u,v,w_1,w_2\}$ and $N_G(v)=\{u,w,v_1,v_2\}$.
By Claim \ref{cl:db-triangle-x}, vertices $v_1,v_2,w_1$ and $w_2$ are all distinct.

\noindent
{\bf Case 1:} $d_G(x)=3$.   Claim \ref{cl:path 3-3-4-3} implies that $d_G(v_1)=d_G(v_2)=d_G(w_1)=d_G(w_2)=4$.
 Denote $N_G(x)=\{u,y,z\}$ and $G'=G-\{v,w,u,x\}$. Then $G'$ has  an  $(h_{G'},M)$-coloring $\psi$.  If $|L(w)\cap L(x)|\geq 2M/31$, then we choose $A\subseteq L(w)\cap L(x)$ such that $|A|=2M/31$. Let $B\subseteq L(x)\setminus A$ such that $|B|=M/31$ as well as $C \subset \psi(w_1) \cup \psi(w_2)$ such that $|C|=M/31$ and $C \cap A=C \cap B=\emptyset$.
By Fact \ref{fact1}, there is a limitation $\psi'$ of $\psi$ on $G'$ such that $B \subset L'(w)$ and $C \subset L'(v)$. We set $C \subset \psi'(v)$, $A \cup B \subset \psi'(w)$, and  $A \cup B \subset \psi'(x)$. Then we add colors to $\psi'(v)$, $\psi'(w)$, and $\psi'(x)$ one by one greedily. Now we can define $\psi'(u)$ properly because each color in $A\cup B$ appears at least twice in neighbors of $u$. A contradiction.

 Let us assume that $|L(w)\cap L(x)|< 2M/31$. Note that $|L(w)| \geq 13M/31$ and $|L(x)| \geq 13M/31$. We assert that $y,z,w_1,$ and $w_2$ are distinct. Otherwise, suppose that $y=w_1$. Then $L(w) \cup L(x) \subset M \setminus \psi(w_1)$ and
$|L(w) \cup L(x)|  \leq 22M/31$. This implies that $|L(w)\cap L(x)| \geq  2M/31$, which is a contradiction and the assertion follows.  Since  $\psi(y)\cup \psi(z) \cup L(x)=M$, we get $|L(w)\cap (\psi(y)\cup \psi(z))|> 11M/31$. This implies that there are $A_1\subseteq L(w)\cap \psi(y)$ and $A_2\subseteq L(w)\cap \psi(z)$ such that $|A_1|=|A_2|=M/31$ and $A_1 \cap A_2=\emptyset$.
 We fix $B \subseteq L(x)$ such that $|B|=M/31$ and $B \cap (A_1 \cup A_2)=\emptyset$. In addition, we choose  $C \subset \psi(w_1) \cup \psi(w_2)$ such that $|C|=M/31$ and $C \cap (A_1 \cup A_2 \cup B)=\emptyset$. By Fact \ref{fact1}, there is a limitation $\psi'$ of $\psi$ on $G'$ such that $ B  \subset L'(w)$ and $C \subset L'(v)$. Deleting colors of $A_1$ first for $\psi'(y)$ and colors of $A_2$ first for $\psi'(z)$, we can assume further $A_1 \cup A_2  \subset L'(x)$.  We put   $A_1 \cup A_2 \cup B \subset \psi'(x)$, $A_1 \cup A_2 \cup B  \subset \psi'(w)$, and $C \subset \psi'(v)$ followed by adding colors to $\psi'(v)$, $\psi'(w)$, and $\psi'(x)$ one by one greedily. As the choices of $A_1,A_2,$ and $B$, we can color $u$ properly to get an $(h_G,M)$-coloring of $G$, a contradiction.

\noindent
{\bf Case 2:} $d_G(x)=4$.  If each vertex from  $\{v_1,v_2,w_1,w_2\}$ have degree four, then let $G'=G-\{u,v,w\}$. Suppose that $\psi$ is  an $(h_{G'},M)$-coloring of $G'$. If $|\psi(x) \cap (\psi(v_1) \cup \psi(v_2))| \geq M/31$, then let $A \subset \psi(x) \cap (\psi(v_1) \cup \psi(v_2))$ with $|A|=M/31$. We choose $B \subset \psi(x)$ such that $A \cap B=\emptyset$ and $|B|=M/31$. By Fact \ref{fact1}, there is a limitation $\psi'$ of $\psi$ on  $G'$ such that $A \subset L'(w)$ and $B \subset L'(v)$. Clearly, $A \cup B \subset \psi'(x)$.  We set  $A \subset \psi'(w)$ and $B \subset \psi'(v)$ followed by adding colors to $\psi'(w)$, $\psi'(v)$, and $\psi'(u)$ one by one greedily.
Then we can define $\psi'(u)$ properly because of  sets $A$ and $B$. We obtain an $(h_G,M)$-coloring of $G$, a contradiction. Thus we need only to consider the case where $|\psi(x) \cap (\psi(v_1) \cup \psi(v_2))| < M/31$. In this case, $|L(v)\cap \psi(x)|>8M/31$ as $L(v) \cup (\psi(v_1) \cup \psi(v_2))=M$.
 We choose $C \subset (\psi(w_1) \cup \psi(w_2))$ such that $|C|=M/31$.
 By Fact \ref{fact1}, there is a limitation $\psi'$ of $\psi$ on  $G'$ such that $C\subseteq L'(v)$.
 We define $\psi'(v)$ and $\psi'(x)$ such that $C \subset \psi'(v)$ and $|\psi'(v) \cap \psi'(x)| \geq 2M/31$.
 To complete the extension of $\psi'$, we pick colors for $w$ and $u$ one by one. Therefore, we can construct an $(h_G,M)$-coloring for $G$, which is a contradiction.

For the case where  one of $v_1,v_2,w_1$, and $w_2$ has degree three, we suppose that $d_G(v_1)=3$.
 Let $v_1'$ and $v_1''$ be the other two neighbors of $v_1$. Then $d_G(v_1')=d_G(v_1'')=4$ by Claim \ref{cl:path 3-3-4-3}.
 Let $G'=G-\{u,v,w,v_1\}$ and $\psi$ be an $(h_{G'},M)$-coloring of $G'$.
 If $|\psi(v_1')\cap \psi(v_1'')|\geq 2M/31$, then we fix $A \subset \psi(v_1')\cap \psi(v_1'')$ with $|A|=2M/31$. Let $B \subset \psi(x)$ and $C \subset \psi(x)$ such that $|B|=|C|=M/31$ as well as $A,B,C$ are pairwise disjoint. By Fact \ref{fact1}, there is a limitation $\psi'$ of $\psi$ on $G-\{u,v,w,v_1\}$ such that $B \subset L'(w)$ and $C \subset L'(v)$. Furthermore, it satisfies that $B \cup C \subset \psi'(x)$ and $A \subset \psi'(v_1') \cap \psi'(v_1'')$. We put $B \subset \psi'(w)$ and $C \subset \psi'(v)$ followed by adding colors to $\psi'(w)$ and $\psi'(v)$ one by one greedily. To complete the extension, we color $v_1$ and $u$ one by one greedily. We obtain an $(h_G,M)$-coloring of $G$, a contradiction.

 If $|\psi(v_1')\cap \psi(v_1'')|< 2M/31$, then we choose $B,C \subset \psi(x)$ as above.
 Additionally, we fix $D \subset (\psi(v_1')\cup \psi(v_1'')) \setminus \psi(v_2)$ such that $|D|=2M/31$ as well as $B,C,$ and $D$ are pairwise disjoint.
 As we did above, we can define a limitation  $\psi'$ of $\psi$ on $G-\{u,v,w,v_1\}$ such that $B \cup C \subset \psi'(x)$, $B \subset L(w)$, and $C \cup D \subset L'(v)$.
We add extra colors to $\psi'(w)$ and $\psi'(v)$ one by one greedily. Then we can get an $(h_G,M)$-coloring for $G$ by picking colors for $v_1$ and $u$ greedily, which is a contradiction. \hfill\rule{1ex}{1ex}


\begin{claim}
\label{cl:db-triangle-3n}
 $G$ contains no $F_4$ such that $u$ has a 3-neighbor.
\end{claim}

\noindent{\bf Proof:}  Claim \ref{cl:triangle-3} implies that each vertex from $\{u,v,w,x\}$ is a 4-vertex. Suppose that $N_G(u)=\{w,v,x,u'\}$ and $N_G(v)=\{w,v,x,v'\}$, where $d_G(u')=3$.  Claim \ref{cl:tr-triangle} gives us that $u' \not =v'$. Moreover, Claim \ref{nop52} implies  that $wu',wv',xu',xv'\notin E(G)$.
Let $G'$ be the graph obtained from $G-\{u,x,w\}$ by identifying $u'$ and $v$.  Claim \ref{cl:c82-v} together with
 Claim \ref{cl:db-triangle-x} yield that
 $G'$ is $C_8^2$-free and $K_4$-free. Therefore, $G'$ has an $(h_{G'},M)$-coloring $\psi$, where we
 assume  $\psi(\underline{u'v})=[1,9M/31]$. We define $A=[1+8M/31,9M/31]$.

 If either $d_G(v')=4$ or $|L(w)\cap L(x)|\geq 2M/31$, then
in the former case, there is a limitation $\psi'$  of $\psi$ on $G-\{u,v,w,x\}$ with $A\subseteq  L'(w) \cap L'(x)$ by Fact \ref{fact1}. Note that $[1,8M/31] \subset L'(v)$. We put $A$ in $\psi'(x)$ and  $\psi'(w)$ followed by adding extra colors to $\psi'(x)$ and  $\psi'(w)$ greedily.
Let $B=[1,8M/31] \setminus (\psi'(x) \cup \psi'(w))$.  We put $B \subset \psi'(v)$ and choose extra colors greedily. This is possible because of  $d_G(v')=4$ and the definition of  $A$. Now, we are able to color $u$ properly because each color from $[1,9M/31]$ appears twice in the neighbors of $u$.
In the latter case, let $\psi'$ be a limitation of $\psi$ on $G-\{u,v,w,x\}$ such that $A\subseteq  L'(w) \cap L'(x)$ and $|L'(w) \cap L'(x)| \geq 2M/31$. To color $x$ and $w$, we  set $A \subset \psi'(x)$ and $A \subset \psi'(w)$ followed by choosing extra colors satisfying $|\psi'(x) \cap \psi'(w)| \geq 2M/31$.
Similarly, let $B=[1,8M/31] \setminus (\psi'(x) \cup \psi'(w))$.  We put $B \subset \psi'(v)$ and pick additional colors for $\psi'(v)$. We are able to construct an $(h_G,M)$-coloring for $G$ by picking colors for $u$ greedily, which is a contradiction.

 If $d_G(v')=3$ and $|L(w)\cap L(x)|<2M/31$,  then let $\psi'$ be a limitation of $\psi$ on $G-\{u,v,w,x\}$ such that $A\subseteq L'(w)\cap L'(x)$. Note that $|(L'(w)\cup L'(x))\cap \psi'(v')| \geq M/31$. We choose $B\subseteq (L'(w)\cup L'(x))\cap \psi'(v')$ with $|B|=M/31$ and denote $B_1=B \cap L'(x)$ as well as $B_2=B \cap L'(w)$.
  We put $A\cup B_1\subset \psi'(x)$ and  $A \cup B_2 \in \psi'(w)$ followed by adding extra colors to $\psi'(x)$ and $\psi'(w)$ greedily. Similarly, to color $v$, we first put $[1,9M/31] \setminus (\psi'(w) \cup \psi'(x))$ in $\psi'(v)$ and choose extra colors greedily.  This is possible because of sets $A$, $B_1$, and $B_2$.   Then we can get an $(h_G,M)$-coloring for $G$ by coloring  $u$  greedily, which is a contradiction.
\hfill\rule{1ex}{1ex}

\begin{claim} \label{noG2}
The graph $G$ contains no $F_7$, see Figure \ref{fig:F7}.
\end{claim}

\noindent
{\bf Proof:}
Suppose that $G$ contains an $F_7$.  Claim \ref{nop52} tells us that $F_7$ is an induced subgraph of $G$.
Furthermore, each vertex in $F_7$ is a 4-vertex by Claim \ref{cl:triangle-3}.  For $z \in \{u,v,a,b\}$, let $z'$ be the neighbor of $z$ not in $F_7$.  Claim \ref{cl:db-triangle-3n} gives us that $u',v',a'$, and $b'$ are 4-vertices. Moreover,
Claim \ref{cl:tr-triangle} implies that $u' \neq v'$ and $a' \neq b'$.

Let $G'$ be the graph obtained from $G-\{x,v,b\}$ by identifying $u$ and $a$. Clearly, Claim \ref{cl:c82-v} and Claim \ref{nop52} tell us that $G'$ is $K_4$-free and $C_8^2$-free. Thus $G'$ has an $(h_{G'},M)$-coloring $\psi$. Suppose that $\psi(\underline{ua})=[1,8M/31]$.
We define  the  limitation $\psi'$  of $\psi$ on $G-\{x,u,v,a,b\}$ with some restrictions. Through removing extra colors properly, we can assume that $|L'(u)|=|L'(a)|=15M/31$, $[1,8M/31] \subset L'(u)$, and $[1,8M/31] \subset L'(a)$.
 Let $A$ be a set  such that $|A|=M/31$ and $A \cap (L'(u) \cup L'(a))=\emptyset$. We can ensure $A \subseteq L'(v)\cap L'(b)$ by deleting colors of $A$ from $\psi(v')$ and $\psi(b')$ first. Similarly, we can assume that $|L'(v)|=|L'(b)|=15M/31$ by removing colors.
We first set   $A \subset \psi'(v)$ and then  define $\psi'(u)$ and $\psi'(v)$ such that $[1,8M/31] \subset  \psi'(u) \cup \psi'(v)$. If $|[1,8M/31] \cap L'(v)|\leq 7M/31$, then we set  $A\cup ([1,8M/31] \cap L'(v))\subset \psi'(v)$ and choose extra colors for $\psi'(v)$ greedily. Note that $|L'(u)\setminus \psi'(v)|\geq 8M/31$ and we can define $\psi'(v)$ properly by including colors from $[1,8M/31] \setminus \psi'(v)$ first. If $|[1,8M/31] \cap L'(v)|>7M/31$, then let $B\subseteq [1,8M/31] \cap L'(v)$ such that $|B|=7M/31$.
Set $\psi'(v)=A\cup B$. To define  $\psi'(u)$, we first put $[1,8M/31]\setminus \psi'(v)$ in $\psi'(u)$ and  add extra colors to $\psi'(u)$ such that $|\psi'(u)|=8M/31$. This is possible because of $|L'(u)\setminus B|=8M/31$.  We define $\psi'(a)$ and $\psi'(b)$ similarly. As each color from $[1,8M/31] \cup A$ appears twice in the neighbors of $x$, we can color $x$ properly.
Therefore, we get an $(h_G,M)$-coloring for $G$, a contradiction. \hfill\rule{1ex}{1ex}

\begin{claim} \label{nocommoneedge}
There are no two triangles sharing an edge.
\end{claim}
\noindent
{\bf Proof:} Suppose that $G$ contains a subgraph consisting of a triangle $uvw$ and a vertex $x$ adjacent to $u$ and $v$.
By Claim \ref{cl:triangle-3}, each vertex from $\{u,v,w,x\}$ is a 4-vertex. Let $u'$ and $v'$ be the other neighbor of $u$ and $v$ respectively. Claim \ref{cl:tr-triangle} tells us that $u' \not =v'$ and
Claim \ref{cl:db-triangle-3n} implies that $d_G(u')=d_G(v')=4$.
Let $G'$ be a graph obtained from $G$ by deleting $u$ and $v$ followed by identifying $w$ and $x$.
 Claim \ref{noG2} together with Claim  \ref{cl:c82-v} yield that $G$ is both $K_4$-free and $C_8^2$-free.
 Thus $G'$ has $(h_{G'},M)$-coloring $\psi$. Note that $|\psi(\underline{wx})|=8M/31$ and $|\psi(u')|=9M/31$. Let $A \subset \psi(u')\setminus \psi(\underline{wx})$ such that $|A|=M/31$. By Fact \ref{fact1}, we may let $\psi'$ be a limitation of $\psi$ on
$G-\{u,v,w,x\}$ such that $A\subseteq L'(v)$. Set $\psi'(w)=\psi(\underline{wx})$ and  $\psi'(x)=\psi(\underline{wx})$. In order to define $\psi'(v)$, we first put $A$ in $\psi'(v)$ and followed by adding extra colors greedily.
 As $d_G(u)=4$ and the definition of $A$, we can color $u$ properly to get an $(h_G,M)$-coloring of $G$, a contradiction.  \hfill\rule{1ex}{1ex}


~

A 4-vertex $v$ is {\it bad} if it is contained in a triangle and has two neighbors with degree three.

\begin{claim} \label{noG6}
$G$ contains no bad 4-vertex.
\end{claim}
\noindent
{\bf Proof:} Suppose  that $G$ contains a bad 4-vertex $v$. Denote $N_G(v)=\{v_1,v_2,v_3,v_4\}$. Suppose that $d_G(v_1)=d_G(v_2)=3$. By Claim \ref{cl:triangle-3}, both $\{v_1,v_2,v_3\}$ and $\{v_1,v_2,v_4\}$ are independent sets. Then we have
$v_3v_4\in E(G)$ and $d_G(v_3)=d_G(v_4)=4$.
For each $1 \leq i \leq 4$, let $v_i'$ and $v_i''$ be the other neighbors of $v_i$. By Claim \ref{cl:path 3-3-4-3} and Claim \ref{cl:triangle-3}, we get $d_G(v_i')=d_G(v_i'')=4$ for $1\leq i \leq 2$.
As Claim \ref{nocommoneedge} shows that there are no two triangles sharing an edge, we get that vertices $v_3',v_3'',v_4',$ and $v_4''$ are distinct.

Let $G'$ be the graph obtained from $G-v$ by identifying $v_1$ and $v_2$. Then  $G'$ is $C_8^2$-free and $K_4$-free by  Claim \ref{cl:c82-v} and  Claim \ref{cl:triangle-3}.
Thus $G'$ has an $(h_{G'},M)$-coloring $\psi$. 
Note that $|L(v_1)|\geq 15M/31$ and $|L(v_2)|\geq 15M/31$. We have $|L(v_1)\cap (\psi(v_3)\cup \psi(v_4))|\geq 2M/31$ and $|L(v_2)\cap (\psi(v_3)\cup \psi(v_4))|\geq 2M/31$. Let  $A\subseteq L(v_1)\cap (\psi(v_3)\cup \psi(v_4))$ such that $|A|=2M/31$. Similarly, let $B\subseteq L(v_2)\cap (\psi(v_3)\cup \psi(v_4))$ such that $|B|=2M/31$.

Let $\psi'$ be a limitation of $\psi$ on $G-\{v,v_1,v_2\}$ such that $A\cup B\subseteq \psi'(v_3)\cup \psi'(v_4)$. We extend $\psi'$ to $v_1$ and $v_2$ such that $A\subseteq \psi'(v_1)$, $B\subseteq \psi'(v_2)$,
and $|\psi'(v_1)\cap \psi'(v_2)|$ as large as possible.  Note that $\psi(\underline{v_1v_2})\subseteq L'(v_1)\cap L'(v_2)$ and $|\psi'(v_1)|=|\psi'(v_2)|=9M/31$. Then $|\psi'(v_1)\cap \psi'(v_2)| \geq 7M/31$.
 Now, $|\cup_{i=1}^4 \psi'(v_i)| \leq 23M/31$ and we can define $\psi'(v)$ properly to obtain an $(h_G,M)$-coloring of $G$, a contradiction.
\hfill\rule{1ex}{1ex}

\begin{claim}\label{noG5}
There are no two triangles in $G$ sharing a vertex.
\end{claim}

\noindent
{\bf Proof:}  Suppose  that $G$ contains two triangles $vv_1v_2$ and $vv_3v_4$ sharing a vertex $v$.
Claim \ref{cl:triangle-3} implies that $d_G(v_i)=4$ for $1 \leq i \leq 4$. We next show that $G$ has an $(h_G,M)$-coloring. The rough idea is that we will start with a coloring of $G-\{v,v_2,v_3\}$, where $v_1$ and $v_4$ are colored with the same colors. We
  extend the coloring to $v_2$, $v_3$, and $v$ one by one. As each $v_i$ has degree four, if there are  $9M/31$ colors such that each of them appears twice in the neighbors of $v$,  then we can color $v$ properly. To achieve this,  when we color $v_3$, we will choose $M/31$ colors from the color set of $v_2$ and assign them to the color set of $v_3$.  In the  process of coloring, we may need to redefine the color sets of $v_1$ and $v_4$.

 For each $1 \leq i \leq 4$,  let $v_i'$ and $v_i''$ be the other two neighbors of $v_i$. By Claim \ref{nocommoneedge}, vertices $v_1',v_1'',v_2'$, and $v_2''$ are distinct. Similarly, $v_3',v_3'',v_4'$, and $v_4''$ are distinct.
 For each $1 \leq i \leq 4$, we have $d_G(v_i')+d_G(v_i'') \geq 7$  by Claim \ref{noG6}.
There are two cases.

\noindent
{\bf Case 1:} There is some $j \in \{1,2\}$ and some $k \in \{3,4\}$ such that $N_G(v_j) \cap N_G(v_k)=v$. By the symmetry, we can assume that $j=1$ and $k=3$. Let $G'$ be the graph obtained from $G-\{v,v_2,v_3\}$ by identifying $v_1$ and $v_4$. Claim \ref{cl:c82-v} and Claim \ref{nocommoneedge} indicate that $G'$ is $C_8^2$-free and  $K_4$-free.
 Thus $G'$ has an $(h_{G'},M)$-coloring $\psi$. Let $\psi'$ be a limitation of $\psi$ on $G-\{v,v_2,v_3,v_2',v_2'',v_3',v_3''\}$.
 Without loss of generality, we assume that $\psi'(v_1)=\psi'(v_4)=[1,8M/31]$.
There are two subcases.

{\bf Subcase 1.1:} $|N_G(v_2) \cap N_G(v_3)|=1$. It implies that $v$ is the only common neighbor of $v_2$ and $v_3$. Notice that $v_2',v_2'',v_3',$ and $v_3''$ are distinct as well as $d_{G'}(u)=d_G(u)-1$ for each $u \in \{v_2',v_2'',v_3',v_3''\}$.
Let $A_2=[1,8M/31] \cap (\psi(v_2') \cup \psi(v_2''))$ and $B_2=\psi(v_2') \cap \psi(v_2'')$. Similarly, let
$A_3=[1,8M/31] \cap (\psi(v_3') \cup \psi(v_3''))$ and $B_3=\psi(v_3') \cap \psi(v_3'')$.
We next color $v_2$, $v_3,$ and $v$ one by one.

Let us start with the case where either $|A_2| \geq  (9-d_G(v_2')-d_G(v_2''))M/31$ or $|B_2| \geq (9-d_G(v_2')-d_G(v_2''))M/31$. In the former case, we delete $M/31$ colors from $\psi(v_2')$ and $M/31$ colors from  $\psi(v_2'')$ respectively, while keeping $A_2$ as large as possible. In the latter case, we delete colors from $\psi(v_2')$ and $\psi(v_2'')$  while keeping $B_2$ as large as possible. In each case,  we can choose colors for $\psi'(v_2)$ greedily. We next color $v_3$. If either $|A_3| \geq  (9-d_G(v_3')-d_G(v_3''))M/31$ or $|B_3| \geq (9-d_G(v_3')-d_G(v_3''))M/31$, then we pick a set $C \subset \psi'(v_2)$ such that $|C|=M/31$. Notice that $C \cap [1,8M/31]=\emptyset$.
 In the latter case, we further require $|C \cap B_3|$ as small as possible.  Let us put $C \subset \psi'(v_3)$. When we define
 $\psi'(v_3')$ and $\psi'(v_3'')$, we
  remove colors from $C$ first while keeping $A_3$ (and $B_3$ in the latter case) as large as possible. In each case, we can choose extra colors for $\psi'(v_3)$ greedily. Note that each color from $[1,8M/31] \cup C$ appears twice in the neighbors of $v$, we can color $v$ properly to get an $(h_G,M)$-coloring of $G$, a contradiction.  Therefore, we  assume that $|A_3|<  (9-d_G(v_3')-d_G(v_3''))M/31$ and $|B_3| < (9-d_G(v_3')-d_G(v_3''))M/31$. Notice that
 $|\psi(v_3')\cup \psi(v_3'')|=(26-d_G(v_3')-d_G(v_3''))M/31-|B_3|.$
 Recall  $d_G(v_4')+d_G(v_4'') \geq 7$. Therefore, $|\psi'(v_4') \cup \psi'(v_4'')| \leq 17M/31$, which yields that there exists a set $D \subset (\psi(v_3') \cup \psi(v_3'')) \setminus (\psi'(v_4') \cup \psi'(v_4''))$ with $|D|=(9-d_G(v_3')-d_G(v_3''))M/31-|B_3|$. The fact $d_G(v_3')+d_G(v_3'') \geq 7$  implies that $|B_3| +|D| \leq 2M/31$. Note that $\psi'(v_2) \cap [1,8M/31]=\emptyset$. Now we are able to choose $C \subset \psi'(v_2)$ such that $|C|=M/31$ and $C \cap B_3=C \cap D =\emptyset$.  Let $E \subset [1,8M/31]$ such that $|E|=|D|$ and $E \cap A_3=E \cap C=\emptyset$. We redefine $\psi'(v_4)$ by removing colors of $E$  and adding colors from $D$. To color $v_3$, we put $C \cup E \subset \psi'(v_3)$.
 We next define $\psi'(v_3')$ and $\psi'(v_3'')$ by removing colors from $C$ first and extra colors if it is needed, here we keep colors from $B_3$. As sets $B_3$ and $D$, we can color $v_3$ properly.  Since each color from $[1,8M/31] \cup C$ appears twice in the neighbors of $v$, we can color $v$ properly to get an $(h_G,M)$-coloring of $G$, a contradiction.

Now, we consider the case where $|A_2| <  (9-d_G(v_2')-d_G(v_2''))M/31$ and $|B_2| < (9-d_G(v_2')-d_G(v_2''))M/31$. In this case, to color $v_2$, we can use the idea to color $v_3$ when $|A_3|<  (9-d_G(v_3')-d_G(v_3''))M/31$ and $|B_3| < (9-d_G(v_3')-d_G(v_3''))M/31$. Namely, we move a set of at most $2M/31$ colors from $[1,8M/31]$ to $\psi'(v_2)$ and replace them by  a set of at most $2M/31$ from $(\psi'(v_2') \cup \psi'(v_2'')) \setminus (\psi'(v_1') \cup \psi'(v_1''))$.
 We repeat the argument above to color $v_3$ by case analysis according to the sizes of $A_3$ and $B_3$.
 No matter whether $|A_3|<  (9-d_G(v_3')-d_G(v_3''))M/31$ or $|B_3| < (9-d_G(v_3')-d_G(v_3''))M/31$, the existence of a desired subset $C \subset \psi'(v_2) \setminus [1,8M/31]$ is ensured by facts   $|D|+|B_3| \leq 2M/31$ and
  $|\psi'(v_2) \setminus [1,8M/31]| \geq 6M/31$.
 As above, each color from $[1,8M/31] \cup C$ appears twice in the neighbors of $v$. Thus we can color $v$ greedily and obtain an $(h_G,M)$-coloring of $G$, a contradiction.

{\bf Subcase 1.2:} $|N_G(v_2) \cap N_G(v_3)| \geq 2$.  We define $A_2,B_2,A_3,$ and $B_3$ as we did in Subcase 1.1. We can reuse ideas in Subcase 1.1 to extend the coloring $\psi'$ to $G$.
For the case where $v_2'=v_3'$ and $v_2'' \not = v_3''$, note that $|\psi(v_2')|=(14-d_G(v_2'))M/31$ as both $v_2$ and $v_3$ are adjacent to $v_2'$. After we finish the coloring of $v_2'$, if we still use $\psi(v_2')$ to denote the color set of $v_2'$, then we require $|\psi(v_2')|=(13-d_G(v_2'))M/31$. Otherwise, we have troubles to extend the coloring to $v_3$. Therefore, to color $v_2$,  we need to consider cases depending on whether $|A_2| \geq (10-d_G(v_2')-d_G(v_2''))M/31$ or $|B_2| \geq  (10-d_G(v_2')-d_G(v_2''))M/31$. Note that when $|A_2| < (10-d_G(v_2')-d_G(v_2''))M/31$ and $|B_2| <  (10-d_G(v_2')-d_G(v_2''))M/31$, it holds $|(\psi(v_2') \cup \psi(v_2'')) \setminus (\psi(v_1') \cup \psi(v_1''))| \geq (10-d_G(v_2')-d_G(v_2''))M/31-|B_2|$ as $|\psi(v_2')|=(14-d_G(v_2'))M/31$ and $|\psi(v_2'')|=(13-d_G(v_2''))M/31$.
 For the case where $v_2'=v_3'$ and $v_2''  = v_3''$,  to define $\psi'(v_2)$, we have to consider cases depending on whether $|A_2| \geq (11-d_G(v_2')-d_G(v_2''))M/31$ or $|B_2| \geq  (11-d_G(v_2')-d_G(v_2''))M/31$ similarly. In each case, we can construct an $(h_G,M)$-coloring of $G$, a contradiction.

\noindent
{\bf Case 2:} For each $j \in \{1,2\}$ and each $k \in \{3,4\}$, we have $|N_G(v_j) \cap N_G(v_k)| \geq 2$. Recall Claim \ref{noG5}. In this case, there are distinct vertices $x,y,w,z$ such that $x \in N_G(v_1) \cap N_G(v_3)$, $y \in N_G(v_1) \cap N_G(v_4)$,
$w \in N_G(v_2) \cap N_G(v_3)$, and $z \in N_G(v_2) \cap N_G(v_4)$. By the symmetry, we assume that $d_G(x)+d_G(y) \geq d_G(w)+d_G(z)$.   Let $G'$ be the graph obtained from $G-\{v,v_2,v_3\}$ by identifying $v_1$ and $v_4$. It is easy to check that $G'$ is $K_4$-free and $C_8^2$-free. Let $\psi$ be an $(h_{G'},M)$-coloring of $G'$ and $\psi'$ be a limitation of $\psi$ on $G-\{v,v_2,v_3,x,w,z\}$. We also assume that $\psi'(v_1)=\psi'(v_4)=[1,8M/31]$ and define $A_2=[1,8M/31] \cap (\psi(w) \cup \psi(z))$ as well as $B_2=\psi(w) \cap \psi(z)$. To color $v_2$, we reuse ideas in Subcase 1.1 and  consider cases depending on whether $|A_2| \geq  (10-d_G(w)-d_G(z))M/31$ or $|B_2| \geq (10-d_G(w)-d_G(z))M/31$. After the coloring of $\psi'(v_2)$, we still use $\psi(w)$ to denote the color set of $w$, here we similarly require   $|\psi(w)|=(13-d_G(w))M/31$. To define $\psi'(v_2')$, when $|A_2| < (10-d_G(w)-d_G(z))M/31$ and $|B_2| < (10-d_G(w)-d_G(z))M/31$, it holds $|(\psi(w) \cup \psi(z)) \setminus (\psi(x) \cup \psi'(y))| \geq (10-d_G(w)-d_G(z))M/31-|B_2|$ as the assumption $d_G(x)+d_G(y) \geq d_G(w)+d_G(z)$.
 To color $v_3$, we define $A_3=[1,8M/31] \cap (\psi(x) \cup \psi(w))$ and $B_3= \psi(x) \cap \psi(w)$ followed by considering cases depending on  $|A_3|$ and $|B_3|$. In each case, we can guarantee that each color of $[1,8M/31] \cup C$ appears twice in the neighbors of $v$ for some $C \subset \psi'(v_2) \setminus [1,8M/31]$ with $|C|=M/31$. Thus we can color $v$ greedily to get an $(h_G,M)$-coloring of $G$, a contradiction.
  \hfill\rule{1ex}{1ex}

\begin{claim} \label{4-4neighbor}
Each 4-vertex has at least one 4-neighbor.
\end{claim}

\noindent{\bf Proof:} Let $v$ be a 4-vertex of $G$ with $N_G(v)=\{v_1,v_2,v_3,v_4\}$. Suppose  that  $d_G(v_i)=3$ for $1 \leq i \leq 4$. Then  $\{v_1,v_2,v_3,v_4\}$ is an independent set and the neighbors of $v_1,v_2,v_3,v_4$ are 4-vertices by Claim \ref{cl:path 3-3-4-3} and Claim \ref{cl:triangle-3}.
Let $G'$ be the graph obtained from $G-v$ by identifying $v_1$ and $v_2$.
 Claim \ref{cl:triangle-3} and  Claim \ref{cl:c82-v} imply that $G'$ is $K_4$-free and contains no $C_8^2$.
Then $G'$ has an $(h_{G'},M)$-coloring $\psi$. Note that $|\psi(\underline{v_1v_2})|=8M/31$. Let $\psi'$ be a limitation of $\psi$ on $G-\{v,v_1,v_2,v_3,v_4\}$. To define $\psi'(v_1)$ and $\psi'(v_2)$, we add  $M/31$ extra colors to $\psi(\underline{v_1v_2})$ greedily.  We can assume that $|L'(v)|=21M/31$, $|L'(v_3)|=15M/31$ and $|L'(v_4)|=15M/31$ by deleting  colors if it is necessary. Let $A\subseteq L'(v)\setminus L'(v_3)$  with $|A|=2M/31$ and $B\subseteq L'(v)\setminus L'(v_4)$ with $|B|=2M/31$.  To color $v$, we put $A\cup B$ in $\psi'(v)$ and add extra colors to
$\psi'(v)$ such that $|\psi'(v)|=8M/31$. Now  $v_3$ and $v_4$  can be colored properly because of sets $A$ and $B$.
We obtained an $(h_G,M)$-coloring of $G$, a contradiction.
%
%
\hfill\rule{1ex}{1ex}

~~~~~

\noindent
{\bf Proof of Lemma \ref{key4}:} Claim \ref{key0} implies that $d_G(v)\geq 3$ for each $v \in V(G)$.
 If $d_G(v)=3$, then   $N_G(v)$ is an independent set by Claim \ref{cl:triangle-3}. Consequently,  $e(N_G(v),N_G^2(v))\geq 6$. If  $d_G(v)=4$, then $v$ has at least a 4-neighbor by Claim \ref{4-4neighbor}. If $v$ has exactly one 4-neighbor, then $v$ is not in any triangles by Claim \ref{cl:triangle-3}, which yields   $e(N_G(v),N_G^2(v))\geq 9$. If $v$ has exactly two 4-neighbors, then $v$ is not in any triangles by Claims \ref{cl:triangle-3} and \ref{noG6}. It follows that  $e(N_G(v),N_G^2(v))\geq 2+2+3+3=10$. If $v$ has at least three 4-neighbors, then $v$ is in at most one triangle by Claim  \ref{nocommoneedge} and Claim \ref{noG5} and we have  $e(N_G(v),N_G^2(v))\geq 9$.  The lemma is proved. \hfill\rule{1ex}{1ex}


\begin{thebibliography}{9}
\bibitem{bk}
O.~Borodin  and A.~Kostochka, On an upper bound on a graph's
chromatic number, depending on the graph's degree and density,
{\it J.~Combin.~Theory Ser.~B,} {\bf 23} (1977), 247--250.
\bibitem{DSV}
 Z.~Dvo\v r\'ak, J.~S.~Sereni, and J.~Volec, Subcubic triangle-free graphs have fractional chromatic number at most $14/5$, {\it J.~Lond.~Math.~Soc.~(2),} {\bf 89}(3) (2014), 641--662.
\bibitem{EK}
K.~Edwards and A.~D.~King, Bounding the fractional chromatic numrber of $K_{\Delta}$-free graphs, {\it SIAM J.~Discrete Math.,} {\bf 27}(2) (2013), 1184--1208.
\bibitem{faj}
S.~Fajtlowicz, On the size of independent sets in graphs, in {\it Proceedings of the Ninth Southeastern  Conference on Combinatorics, Graph Theory and Computing}(Florida  Atlantic Univ, Boca Raton, FL,1978), Congress, Numer, XXI, Utilitas  Math, Winnipeg, MB, 1978, 269--274.
\bibitem{FKK}
D.~Ferguson, T.~Kaiser, and D.~Kr\'al',  The fractional chromatic number of triangle-free subcubic graphs, {\it European J.~Combin.,} {\bf 35} (2014), 184--220.
\bibitem{gallai}
T.~Gallai, Kritische graphen I, {\it Magyar Tud. Akad. Mat. Kutat \'o Int. K\"ozl}, 8 (1963), 165--192.
\bibitem{gm}
J.~Griggs and O.~Murphy, Edge density and independence ratio in triangle-free graphs with maximum degree three, {\it Discrete Math.}, {\bf 152} (1996), 157--170.
\bibitem{hz}
H.~Hatami and X.~Zhu, The fractional chromatic number of graphs of maximum degree at most three, {\it SIAM J.~Discrete Math.,} {\bf 24}(4) (2009/2010), 1762--1775.
\bibitem {ht}
C.~C.~Heckman and R.~Thomas, A new proof of the independence ratio of triangle-free cubic graphs, {\it Discrete Math.}, {\bf 233} (2001), 233--237.
\bibitem{johansson}
A.~Johansson, Asymptotic choice number for triangle free graphs, {\it DIMACS Technical Report,} {\bf 91-5}.
\bibitem{KKV}
F.~Kardo\v s, D.~Kr\'al', and J.~Volec,  Fractional colorings of cubic graphs with large girth, {\it SIAM J.~Discrete Math.,} {\bf 25}(3) (2011), 1454--1476.
\bibitem{KLP}
 A.~D.~King, L.~Lu, and X.~Peng, A fractional analogue of Brooks' Theorem,  {\it SIAM.~J.~Discrete Math.,} {\bf 26}-2 (2012), 452--471.
\bibitem{chl}
C.~H.~Liu, An upper bound on the fractional chromatic number of triangle-free subcubic graphs. {\it SIAM J.~Discrete Math.,} {\bf 28}(3) (2014), 1102--1136.
\bibitem{LP}
L.~Lu and X.~Peng, The fractional chromatic number of triangle-freegraphs with $\Delta \leq 3$, {\it Discrete Math.,} {\bf 312}(24) (2012), 3502--3516.
\bibitem{mr}
M.~Molloy and B.~Reed, {\it Graph colouring and the probabilistic method},
volume 23 of Algorithms and Combinatorics. Springer-Verlag, Berlin, 2002.
\bibitem{reed}
B.~Reed, $\omega, \Delta,$ and $\chi$, {\it J.~of Graph Theory}, vol {\bf 27}-4 (1998), 177--227.
\bibitem{reed1}
B.~Reed, A strengthening of Brooks' thoerem, {\it J.~Combin.~Theory Ser.~B,}  {\bf 76} (1999), 136--149.
\bibitem{su}
E.~R.~Scheinerman and D.~H.~Ullman, {\it Fractional Graph Theory}. A Rational Approach to the Theory of Graphs,  Wiley-Intersci. Ser. Discrete Math. Optim, John Wiley \& Sons, Inc,  New York, 1997.
\bibitem{staton}
W.~Staton, Some Ramsey-type numbers and the  independence ration, {\it Trans.~Amer.~Math.~Soc.}, {\bf 256} (1979), 353--370.




\end{thebibliography}
\end{document}